\documentclass[%
 aip,
 amsmath,amssymb,
 reprint,%
]{revtex4-1}

\usepackage{graphicx}% Include figure files
\usepackage{dcolumn}% Align table columns on decimal point
\usepackage{bm}% bold math
\usepackage[utf8]{inputenc}
\usepackage[T1]{fontenc}
\usepackage{mathptmx}
\usepackage{color}

\begin{document}

\preprint{AIP/123-QED}

\title[Leonid Shilnikov and mathematical theory of dynamical chaos]{Leonid Shilnikov and mathematical theory of dynamical chaos}

\author{Sergey Gonchenko}
\email{sergey.gonchenko@mail.ru}
\affiliation{
Scientific and Educational Mathematical Center ``Mathematics of Future Technologies'',
Lobachevsky State University of Nizhny Novgorod, \\
23 Gagarina Ave., 603950 Nizhny Novgorod, Russia}
\author{Alexey Kazakov}
\email{kazakovdz@yandex.ru}
\affiliation{National Research University Higher School of Economics, \\
25/12 Bolshaya Pecherskaya Ulitsa, 603155 Nizhny Novgorod, Russia}
\author{Dmitry Turaev}
\email{d.turaev@imperial.ac.uk}
\affiliation{Imperial College, London, UK}
\affiliation{National Research University Higher School of Economics, \\
25/12 Bolshaya Pecherskaya Ulitsa, 603155 Nizhny Novgorod, Russia}

\author{Andrey L. Shilnikov}
\email{ashilnikov@gsu.edu}
\affiliation{Neuroscience Institute and Department of Mathematics \& Statistics, Georgia State University, \\
100 Piedmont Ave., Atlanta, GA 30303, USA.}
%\affiliation{National Research University Higher School of Economics, \\
%25/12 Bolshaya Pecherskaya Ulitsa, 603155 Nizhny Novgorod, Russia}

\date{\today}% It is always \today, today,

%\begin{abstract}
%\end{abstract}

\maketitle

\begin{figure}[t!]
\includegraphics[width=.9\linewidth]{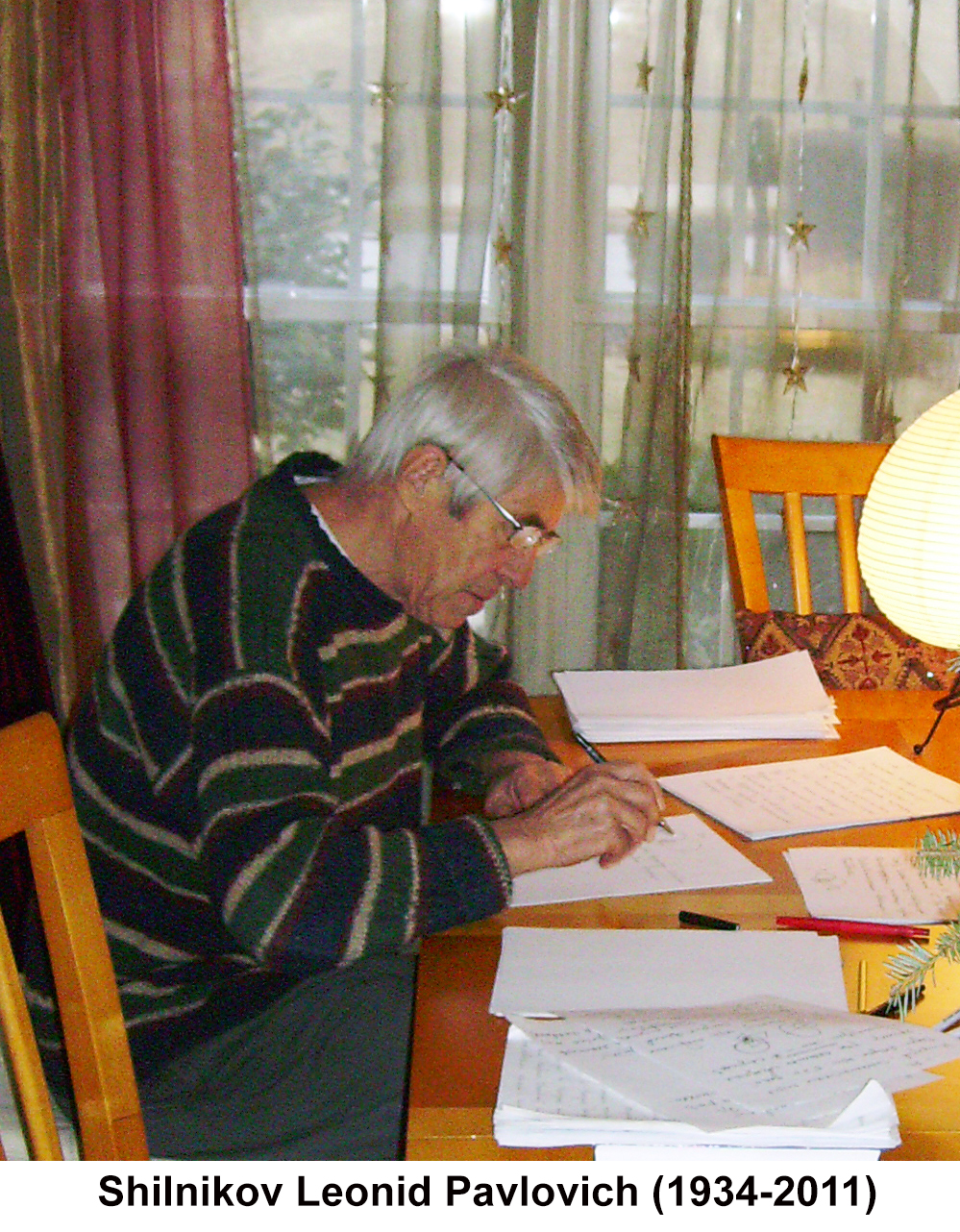}
\end{figure}

\section*{Introduction}

This Focus Issue ``Global Bifurcations, Chaos, and Hyperchaos: Theory and Applications'' is dedicated to the 85th anniversary of the great mathematician, one of the founding fathers of dynamical chaos theory, Leonid Pavlovich Shilnikov. His works decisively influenced the nonlinear dynamics. Research of many scholars from various sciences is deeply rooted in Shilnikov's legacy.

In 1965, he discovered and described what is now known as spiral, or Shilnikov chaos, and since that he had made the study of global bifurcations and transition to chaos the central topic of his work.
Shilnikov, and the dedicated scientific group he had created, investigated a multitude of homoclinic phenomena, gave the first comprehensive theoretical description of the classical Lorenz attractor, discovered the synchronization-to-chaos transition via the torus breakdown. Shilnikov's pioneering results have long become classical and are included in text- and reference books on dynamical systems and bifurcations.

As a true great scientist, L.P.~Shilnikov had a ``magical vision'' gift that let him establish connections between seemingly unrelated topics, find an unexpected formulation of a problem, or propose an approach that would suddenly become evident only many years later. Perhaps, because of that he became a global ``attractor'' for many colleagues and researchers from mathematics to physics, biology, neuroscience, chemistry, and engineering, who deeply valued personal and professional contacts with L.P.~Shilnikov. Many scholars acknowledge that his ideas and charisma vastly influenced their own development.

For all of us, Leonid Pavlovich was more than just a wise advisor. He was a wonderful individual, our dear colleague and a role model in science and life. We miss him greatly. As a tribute to his memory, we put together this Focus Issue with the overall idea to highlight a recent progress in the study of higher-dimensional dynamics done by Shilnikov's students and followers.  This goes back to our many conversations with Shilnikov. He stressed that chaotic dynamics were reasonably understood mostly for three-dimensional systems of differential equations or two-dimensional maps and that the focus of the dynamical systems research must move to higher dimensions.

We decided to start with an overview of Shilnikov's life work on the theory of dynamical systems, which will be followed by the description of the papers contributed to this Focus Issue. We selected to discuss the following themes:
\begin{itemize}\itemsep 1pt
\item Global bifurcations of high-dimensional dynamical systems,
\item Shilnikov saddle-focus and spiral chaos,
\item Homoclinic chaos,
\item Homoclinic tangency,
\item Mathematical theory of synchronization,
\item Lorenz attractor, quasiattractors, and pseudohyperbolicity.
\end{itemize}
His interests were much broader and included a variety of other topics, such as the theory of Hamiltonian chaos, non-autonomous systems, Markov chains and symbolic dynamics of hyperbolic systems, the theory of solitary waves and its applications, to name a few. The part of Shilnikov legacy, which we have chosen to elucidate, contains the most pioneering and fundamental results which made the greatest impact on the chaos theory and determined the major directions for the future research.

\section{Early works: Global bifurcations}
\label{sec:early}
The theory of bifurcations had grown from the pioneering works of A.A.~Andronov and E.A.~Leontovich in the 30s of the XX-th century \cite{andronov1937some, andronov1938, TB1, TB2}. Motivated by the development of radiophysics (the theory of oscillations), they showed that a stable periodic orbit of a system of ordinary differential equations on a plane can be born either out of a semi-stable periodic orbit, or from an equilibrium state in what is now called the Andronov-Hopf bifurcation, or from a homoclinic loop to a saddle or to a saddle-node. Their results relied on the previously developed Poincar\'e-Bendixson-Dulac theory of systems on the plane. The mathematical methods for the analysis of multi-dimensional systems were not developed at the time, so the higher-dimensional theory of bifurcations was reserved for the next generations.

Some 20 years later, Shilnikov became interested in the problem of the generalization of the Andronov-Leontovich theory of homoclinic loops onto systems of an arbitrary dimension. The challenge was to develop new theoretical approaches and mathematical techniques to generalize the planar results for high-dimensional systems; no one expected the groundbreaking discoveries which he was going to make.  

\begin{figure}[t!]
\includegraphics[width=.99\linewidth]{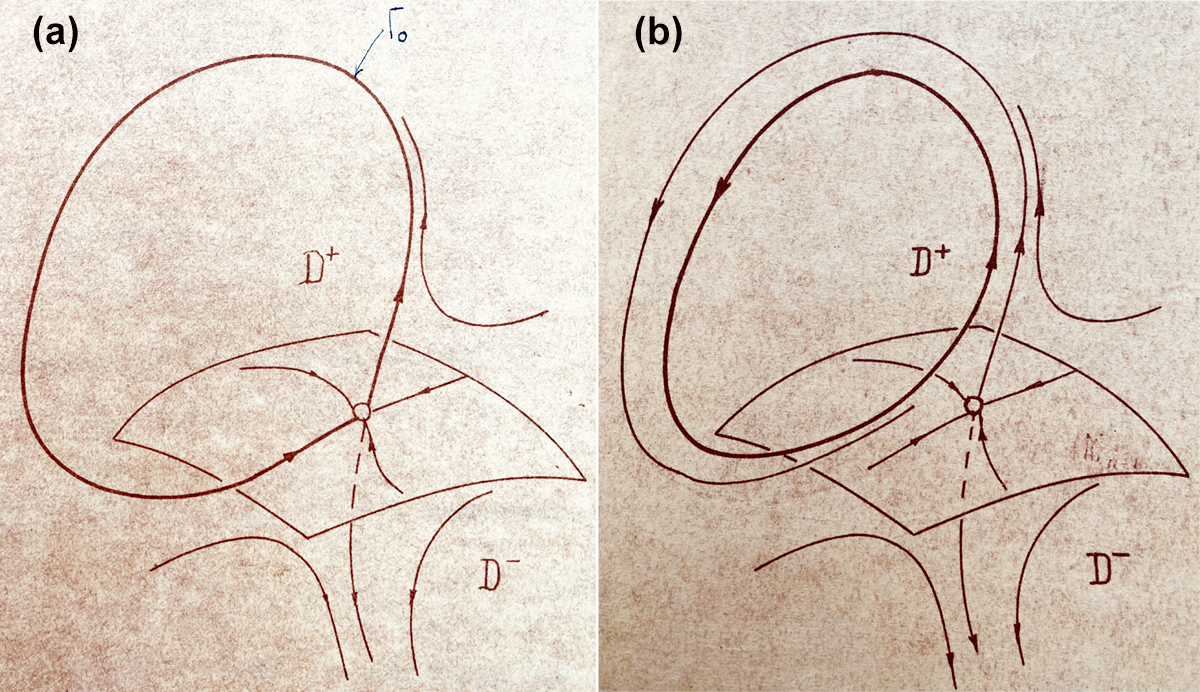}
\caption{Bifurcation of a homoclinic orbit $\Gamma_0$ to a saddle giving rise to the emergence of a single stable periodic orbit provided that the characteristic exponents fulfill the condition $-\mbox{{\rm Re}}\;\lambda_i(0)>\lambda_n(0), \; i = 1,...,n - 1$. This and the next figure are hand-drawings from Shilnikov Ph.D. thesis \cite{Shil_thesis}. }
\label{fig1}
\end{figure}

\begin{figure}[t!]
\includegraphics[width=.99\linewidth]{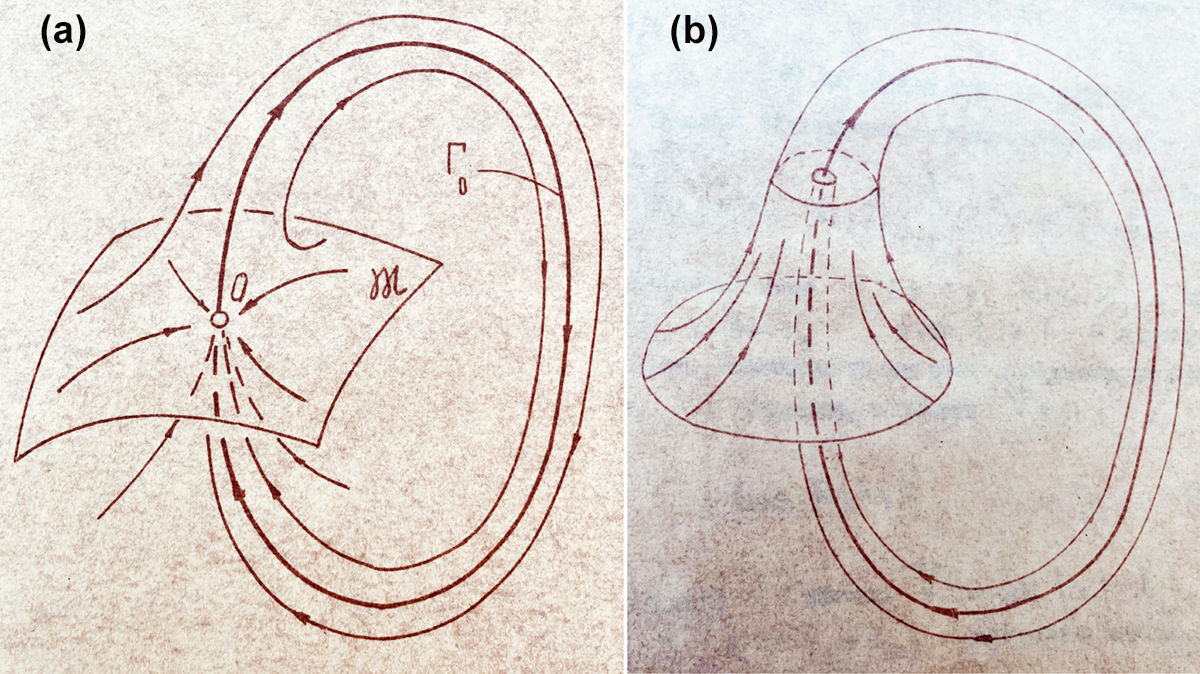}
\caption{Bifurcation of a homoclinic saddle-node equilibrium state resulting in the appearance of a single stable periodic orbit after the saddle-node is gone.}
\label{fig2}
\end{figure}

First, he studied two homoclinic bifurcations that result in the emergence of a single stable periodic orbit out of a homoclinic loop. These are:

\begin{enumerate}
\item The bifurcation of a homoclinic loop $\Gamma_0$ to a saddle equilibrium $O$ with eigenvalues such that $\mbox{Re}\;\lambda_i <0$ $\;(i=1,...,n-1), \lambda_n>0$ and the so-called saddle value $\sigma \equiv \lambda_n + \max \mbox{Re}\;\lambda_i$ is negative; see Fig.~\ref{fig1}.
\item The bifurcation of a homoclinic loop $\Gamma_0$ to a saddle-node equilibrium $O$ satisfying $\mbox{Re} \; \lambda_i <0, i = 1, ..., n-1, \lambda_n = 0$; see Fig.~\ref{fig2}.
\end{enumerate}

These results were published in the papers\cite{Sh62,Sh63} which became chapters of his PhD\cite{Shil_thesis} that was defended in 1962. The results were indeed in line with what happens on the plane. The reason is the strong dissipation (area-contraction property) of the phase space of system near the homoclinic loop. However, the new techniques were ready to be applied to the general case, which Shilnikov did next. The crucial finding was that when the strong dissipation is violated, the behavior near the homoclinic loop depends on whether the eigenvalues of the linearization matrix at the saddle are real or complex, and in the complex case the dynamics near the loop can be chaotic!

This result was published in 1965 \cite{Sh65}; the word ``chaos'' was not a scientific term at that time, so the result was formulated as the existence of infinitely many saddle periodic orbits near the homoclinic loop. This fact, that a simple orbit (a homoclinic to a saddle-focus) can cause complex dynamics, was an amazing and disturbing discovery, one of the critical events (alike the Smale horseshoe and the Anosov torus) that initiated the dynamical chaos theory. We discuss the type of chaotic behavior associated with the Shilnikov saddle-focus loop in the next Section.

\begin{figure}[t!]
\includegraphics[width=.99\linewidth]{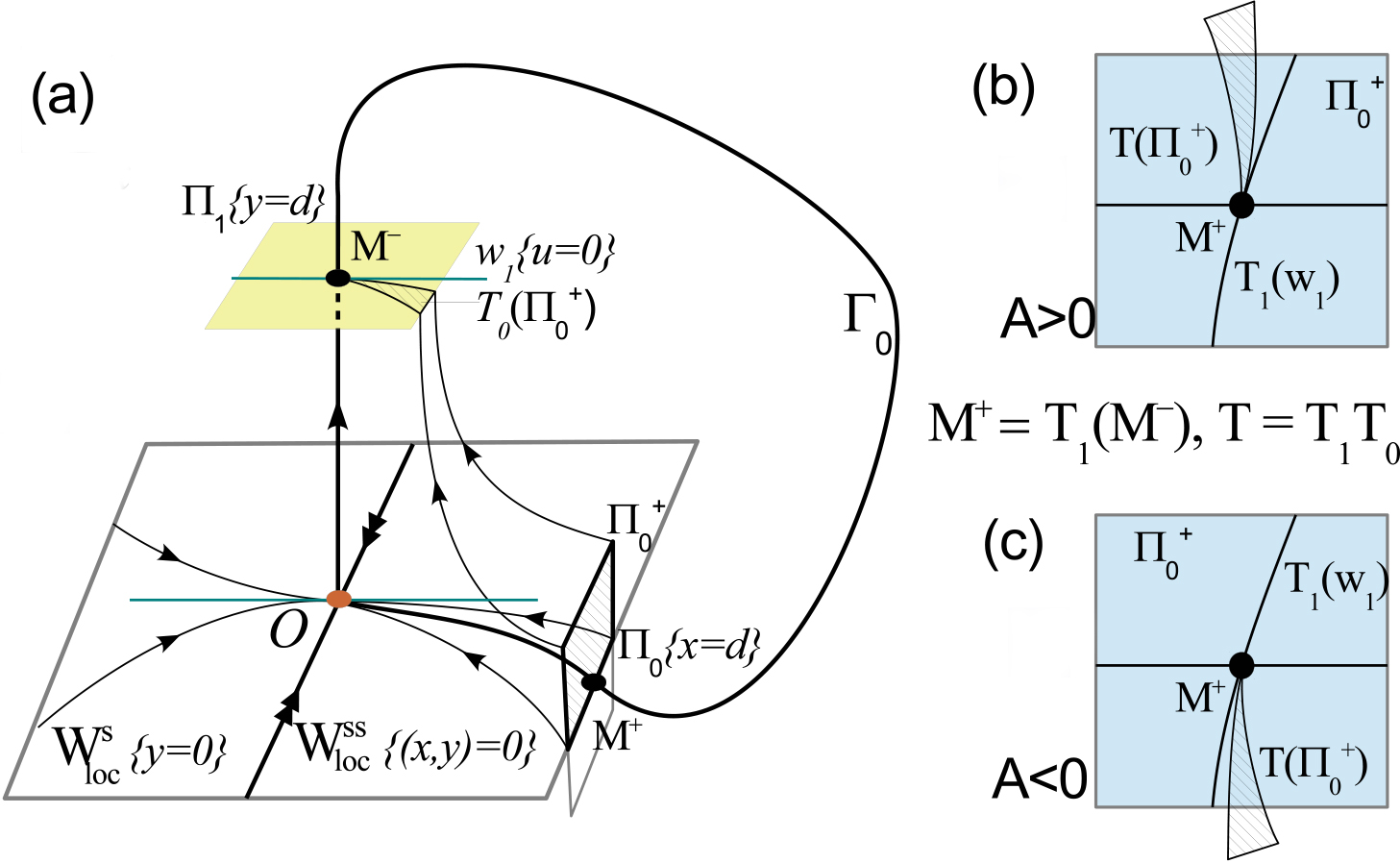}
\caption{(a) Example of a homoclinic loop $\Gamma_0$ to a saddle $O$ with real eigenvalues $\lambda_{ss} <\lambda_{s} <0 <\lambda_u$ in 3-dimensional phase space. If the saddle value $\sigma = \lambda_u + \lambda_s <0$, a stable periodic orbit is born from such loop, see Fig.~\ref{fig1}. When $\sigma$ is positive, then additional conditions are needed for a single saddle periodic orbit to be born. In some local coordinates $(x,y,u)$ near $O$ the unstable and stable manifolds are the $y$-axis and the $(x,u)$-plane; the leading eigen-direction is the $x$-axis, and the eigen-direction corresponding to the non-leading eigenvalue $\lambda^{ss}$ is the $u$-axis. We require that the separatrix $\Gamma_0$ does not lie in the strong-stable manifold $W^{ss}$. Hence, $\Gamma_0$ enters $O$ as $t \to + \infty$ along the $x$-axis. In order to see that the separatrix value $A$ does not vanish, we build two cross-sections to the homoclinic loop: $\Pi_0$ across the local stable manifold $W^s_{loc}$ and $\Pi_1$ across the local unstable manifold. Orbits starting on $\Pi_0$ pass near $O$ and hit $\Pi_1$ at the points of a thin wedge tangent to the $x$-axis; the flow along $\Gamma_0$ returns this wedge to $\Pi_0$. The condition $A\neq 0$ means that the image wedge on $\Pi_0$ is not tangent to $W^s_{loc}$, as shown in panels (b,c) that represent the cases of the orientable (b) and non-orientable (c) homoclinic loop.}
\label{fig3}
\end{figure}

In the other case, when the nearest to the imaginary axis eigenvalue is real, Shilnikov showed\cite{Sh68a} that a single periodic orbit emerges out of a homoclinic loop to a saddle equilibrium provided the four genericity conditions are satisfied: (i) the saddle value $\sigma$ must be non-zero, (ii) the separatrix value $A$ must be non-zero, (iii) when the homoclinic loop enters the saddle it must be tangent to the leading direction, and (iv) there must be only one, and simple, leading eigenvalue, see Fig.~\ref{fig3}.

The papers\cite{Sh65,Sh66,Sh68a} gave an initial impulse to a rich theory of homoclinic phenomena. This theory became an endless source of examples and conceptual models of chaotic behaviors. Many of the results are now a folklore, however they have a single origin -- the entire theory has grown out of the early Shilnikov's works.

\begin{figure}[t!]
\includegraphics[width=.99\linewidth]{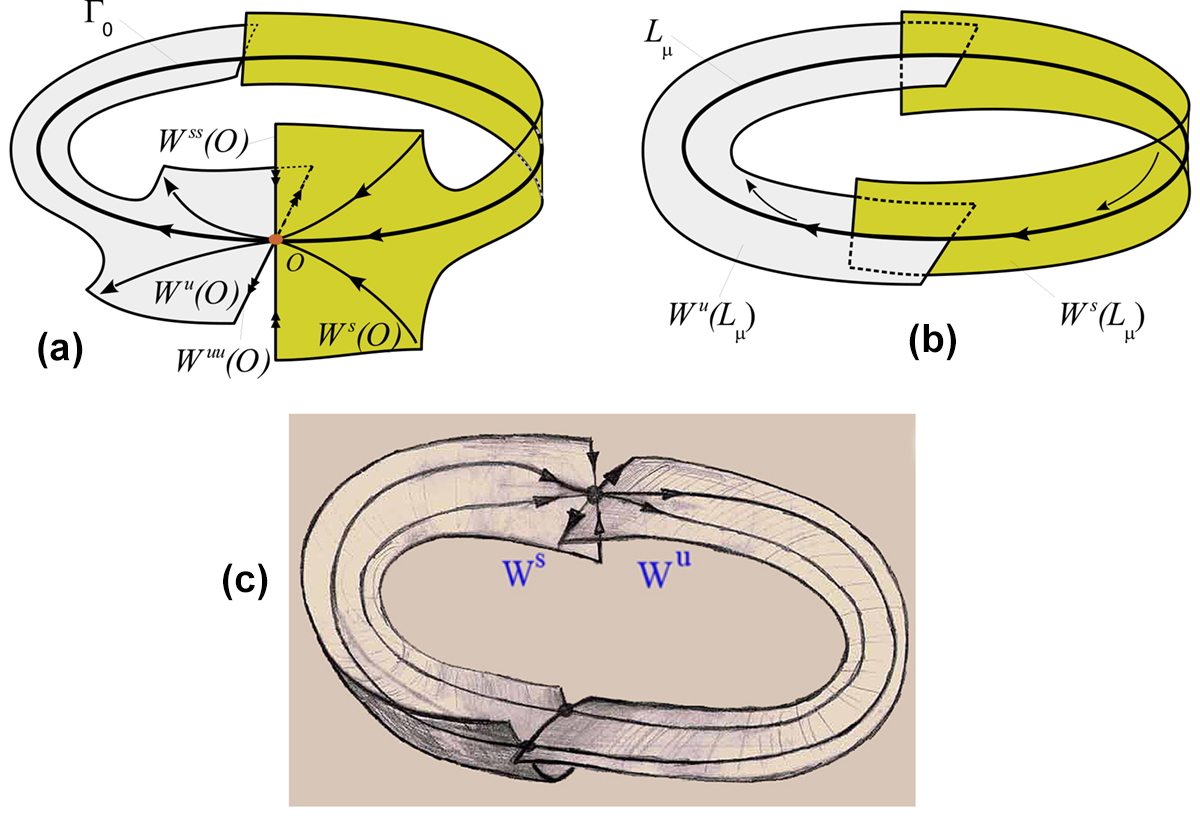}
\caption{(a) A saddle-saddle equilibrium $O$ in 3-dimensional phase space with a single homoclinic orbit $\Gamma_0$, along which the 2-dimensional stable and unstable invariant manifolds $W^s(O)$ and $W^u(O)$ intersect transversely. (b) Once the saddle-saddle vanishes, the homoclinic transitions to a single saddle periodic orbit $L_\mu$. (c) A saddle-saddle can have several transverse homoclinic orbits; from Ref.~[\onlinecite{saddle-saddle}].}
\label{fig4}
\end{figure}

Shilnikov's work on the general case of the homoclinic to a saddle-node is less known to the applied research community, but it also led to a remarkable discovery. In Refs.~\cite{Sh66,Sh69}, he examined the homoclinic bifurcation of a saddle-saddle equilibrium state, see Fig.~\ref{fig4}a. Unlike the homoclinic saddle-node from his early publication\cite{Sh63}, the saddle-saddle has positive real parts in addition to the zero eigenvalue and the eigenvalues with negative real parts. So, both its unstable and stable manifold have dimension larger than one. Shilnikov proved that if (i) the homoclinic loop $\Gamma_0$ leaves the saddle-saddle $O$ and comes back along the eigen-direction corresponding to the zero eigenvalue and (ii) the stable and unstable manifolds of $O$ intersect along $\Gamma_0$ transversely, then a single saddle periodic orbit emerges from the homoclinic orbit to a saddle-saddle after the latter vanishes, see Fig.~\ref{fig4}b. This looks like the saddle-node case, where a stable periodic orbit is now replaced with a saddle one.

\begin{figure}[t!]
\includegraphics[width=.99\linewidth]{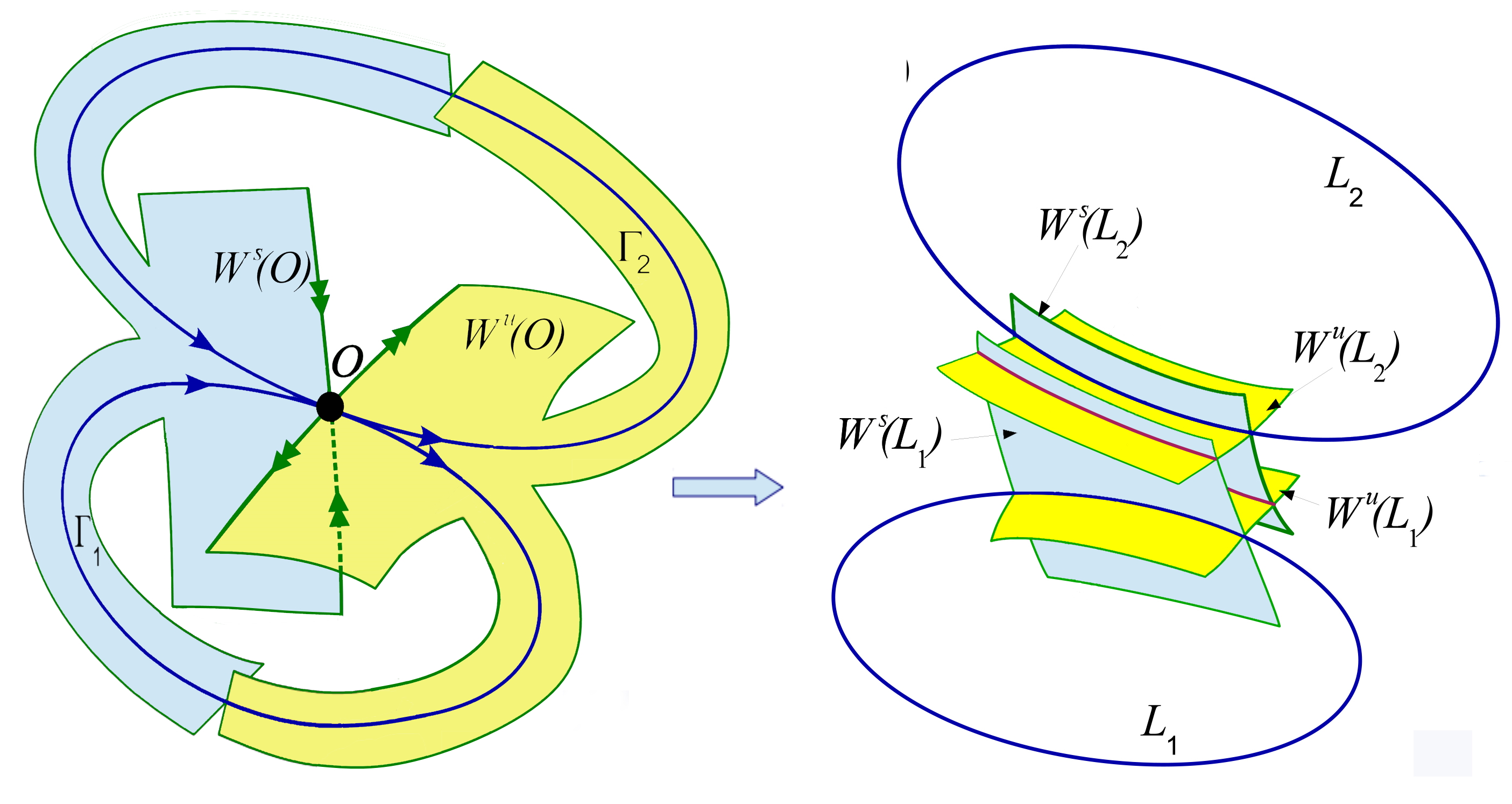}
\caption{As a saddle-saddle $O$ with two homoclinic loops $\Gamma_1$ and $\Gamma_2$ disappears, the stable and unstable manifolds of the saddle cycles $L_1$ and $L_2$ intersect transversely, which causes the onset of a hyperbolic set by virtue of Shilnikov theorem\cite{Sh67Msb}, see Section \ref{sec_3}.}
\label{fig7}
\end{figure}

Let us emphasize that the saddle-saddle can have more than one homoclinic loop, see Fig.~\ref{fig4}c. Shilnikov\cite{Sh69} showed that after a saddle-saddle with $p$ homoclinic loops disappears, then (i) each loop  produces a saddle periodic orbit, and (ii) there emerges {\em a non-trivial hyperbolic set, which is in one-to-one correspondence with orbits of the topological Bernoulli shift on $p$ symbols}, see Fig. \ref{fig7}.

No matter how many homoclinic loops the saddle-saddle has, this bifurcation remains of a codimension-1. Thus, in any one-parameter family which crosses the bifurcation surface corresponding to the Shilnikov saddle-saddle, one should observe a sudden transition from simple dynamics (a pair of saddle equilibria) to chaos (the hyperbolic set after the saddle-saddle is gone).

This was the first example of a \textit{scenario of the transition to chaos}. The search for the transition-to-chaos scenarios and their analysis became one of the main topics of Shilnikov's research for many years: such scenarios explain how chaotic dynamics can arise in real-world applications and determine the further development of the theory. Other scenarios, like period-doublings and torus breakdown, were found later. However, the Shilnikov saddle-saddle scenario is special because the chaotic dynamics is fully formed in just a single bifurcation. Another scenario with the same property was discovered much later -- the birth of a hyperbolic attractor through the blue-sky catastrophe \cite{TSh95}, see Section \ref{sec:tor}.

\section{Spiral chaos}
\label{spiral}
\begin{figure}[t!]
\includegraphics[width=.99\linewidth]{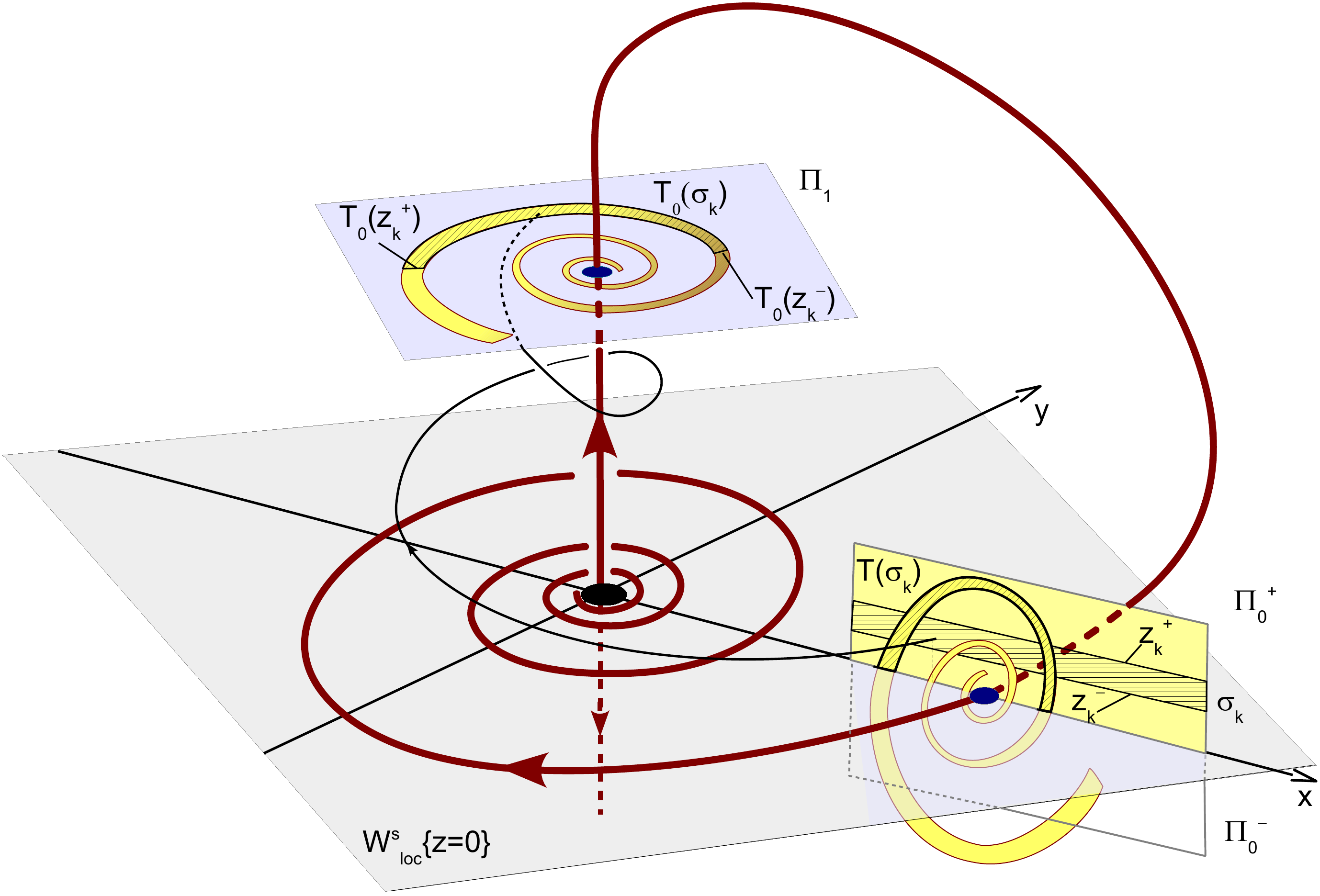}
\caption{Shilnikov saddle-focus loop in 3D phase space: a homoclinic orbit $\Gamma_0$ to a saddle-focus equilibrium $O$ with a two-dimensional stable manifold and a one-dimensional unstable manifold. When the saddle value is positive, there are infinitely many Smale horseshoes in the return map on a cross-section $\Pi^+$ transverse to the local stable manifold $W^s_{loc}(O)$.}
\label{fig5}
\end{figure}

The {\em Shilnikov loop} is a homoclinic orbit to a saddle-focus with positive saddle value. In the pioneering paper\cite{Sh65} Shilnikov showed that the existence of such homoclinic loop in the three-dimensional phase space implies the existence of infinitely many saddle periodic orbits near it. Next, he extended this result onto the four-dimensional case\cite{Sh67sf} where all eigenvalues at the saddle-focus are complex. In the following paper\cite{Sh70} he analyzed the general multi-dimensional case and proved that a neighborhood of the homoclinic loop contains a hyperbolic (chaotic) invariant set which includes infinitely many Smale horseshoes. He also gave a detailed description of the chaotic set and showed that its structure is far more complex than that of the Smale horseshoes. In particular, along with the chaotic set, a neighborhood of the Shilnikov saddle-focus loop can contain infinitely many {\em stable} periodic orbits \cite{ovsyannikov1986systems,ovsyannikov1991systems}.

Figure~\ref{fig5} illustrates the geometric idea behind the proof of the existence of chaos near the saddle-focus loop. We consider the three-dimensional case where the saddle-focus has eigenvalues $\lambda \pm i \omega$ and $ \gamma $ such that $\omega \neq 0$ and $\lambda <0 <\gamma$, so the stable manifold of the saddle-focus is two-dimensional, and the unstable manifold is one-dimensional. Let $\lambda +\gamma > 0$. Here, as in the case of a homoclinic loop to a saddle (see Fig.~\ref{fig3}), the dynamics are determined by the first-return map $T$ of a two-dimensional cross-section $\Pi_0^+$ transverse to the loop $\Gamma_0$. Due to the complex eigenvalues, the image $T(\Pi_0^+)$ has a spiraling shape, see Fig.~\ref{fig5}, so the pre-image $T^{-1}(\Pi_0^+)\cap \Pi_0^+$ is a sequence of disjoint strips $\sigma_k$ accumulating on $W^s_{loc}(O)\cap \Pi_0^+$. The positivity of the saddle value $\lambda+\gamma$ implies that the snake $T(\Pi_0^+)$ is large enough, so that for each $k$ the image $T(\sigma_k)$ intersects $\sigma_k$ and $T$ acts as the Smale horseshoe map on $\sigma_k$. Moreover, $T(\sigma_k)$ intersects other strips too. As Shilnikov proved\cite{Sh70}, this creates an infinite Markov chain whose structure is controlled by the value of $\rho=\lambda/\gamma$. The latter result became a model for the theory of homoclinic tangencies. It also directly implies the sensitive parameter dependence of the fine structure of the Shilnikov chaos, which, in fact, results in the utmost complexity of its dynamics and bifurcations (see Section~\ref{sec_3A}).

The stunning discovery of chaos near a homoclinic to a saddle-focus was a defining moment in Shilnikov's scientific career and the starting point of his life-long quest for homoclinic structures underlying the many faces of dynamical chaos. His works gained a recognition in the West after a series of papers by Arneodo, Coullet, and Tresser \cite{arneodo1980occurence, arneodo1981possible, arneodo1982oscillators} who emphasized the importance of the Shilnikov homoclinic loop for the chaos theory. Starting with late 70s, the spiral chaos -- strange attractors due to the Shilnikov saddle-focus -- keeps been discovered as the most prominent dynamical regime in various models from hydrodynamics, electronics, optics, astrophysics, chemistry, biology, mechanics, neuroscience, etc.\cite{Sciheritage,GGKKB19}.

\begin{figure}[h!]
\includegraphics[width=.5\linewidth]{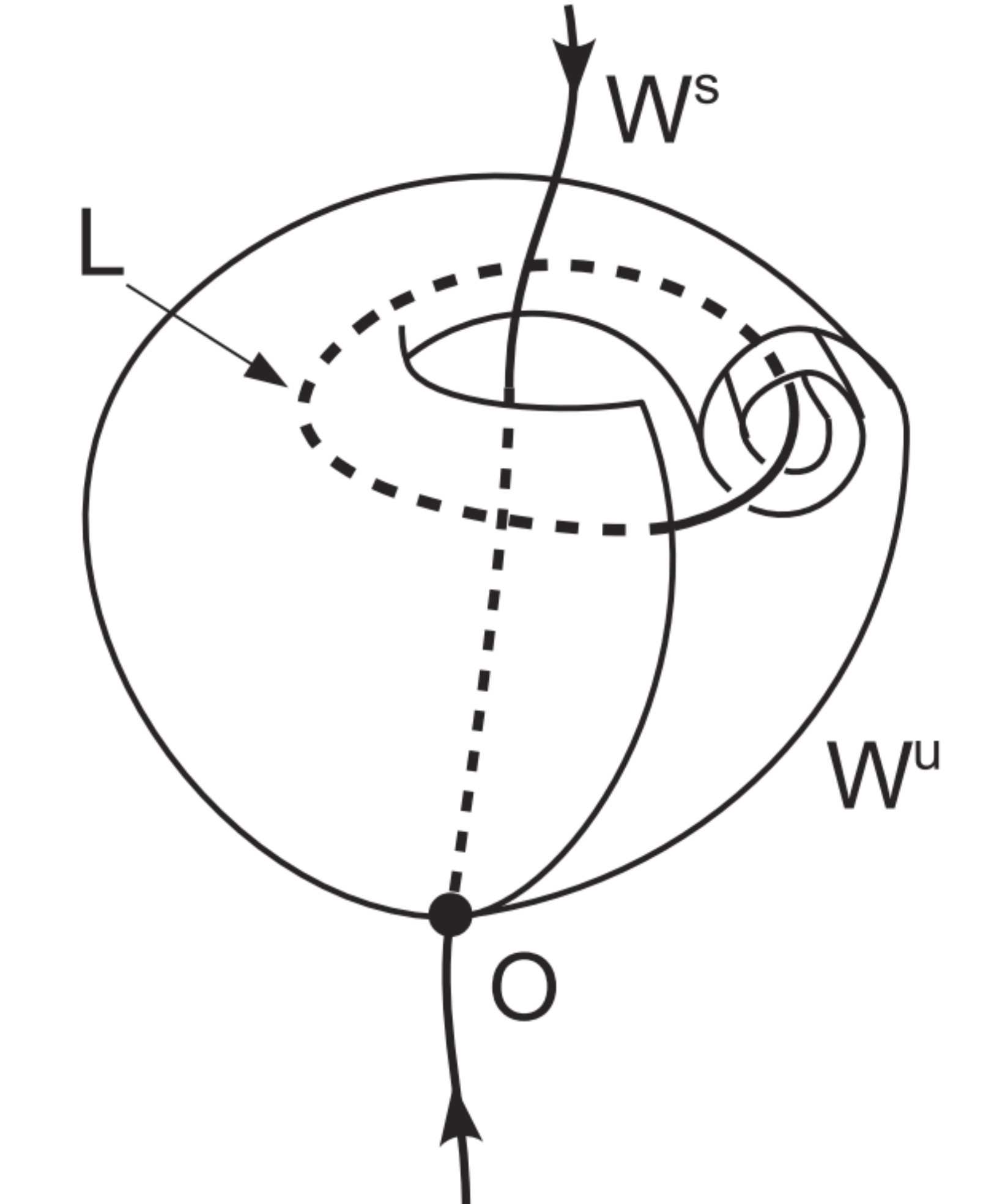}
\caption{Formation of the Shilnikov whirlpool near a saddle-focus equilibrium $O$ with the two-dimensional unstable manifold $W^u$ precedes the homoclinic bifurcation.}
\label{fig65}
\end{figure}

This explosive development got Shilnikov by surprise: why the spiral chaos occurs in so many diverse models, regardless of their nature? As an answer, he proposed {\em a universal scenario of spiral chaos formation}\cite{Sh86} along the following lines. It is common that the system, at some parameter values, operates at a stationary regime (a globally stable equilibrium), but a change in parameters makes the dynamics chaotic. The transition to chaos requires the loss of stability of the equilibrium state, and this usually happens via the Andronov-Hopf bifurcation where a pair of complex-conjugate eigenvalues of the linearized system crosses the imaginary axis. Shilnikov emphasized that the equilibrium becomes a saddle-focus after this bifurcation. As parameters change further, the two-dimensional unstable manifold of the saddle-focus can form the so-called Shilnikov whirlpool, see Fig.~\ref{fig65}, which bounds the basin of an attractor. After that, homoclinic loops can form that give rise to the spiral chaos \cite{Sh86}.

In the same paper\cite{Sh86} Shilnikov described several realizations of this scenario: ``safe'', related to a period-doubling cascade or to a torus breakdown, and also ``dangerous'' ones, where the spiral chaos emerges immediately after the subcritical Andronov-Hopf bifurcations. He proposed several geometric models for them, and detailed phenomenological principles of the onset of other types of chaotic attractors, as was further developed in \cite{GGS12,GGKT14}. In an unpublished manuscript, he built a theory of the emergence of hyperbolic attractors through saddle-focus homoclinic bifurcations.

\section{Homoclinic chaos} \label{sec_3}

\begin{figure}[h!]
\includegraphics[width=\linewidth]{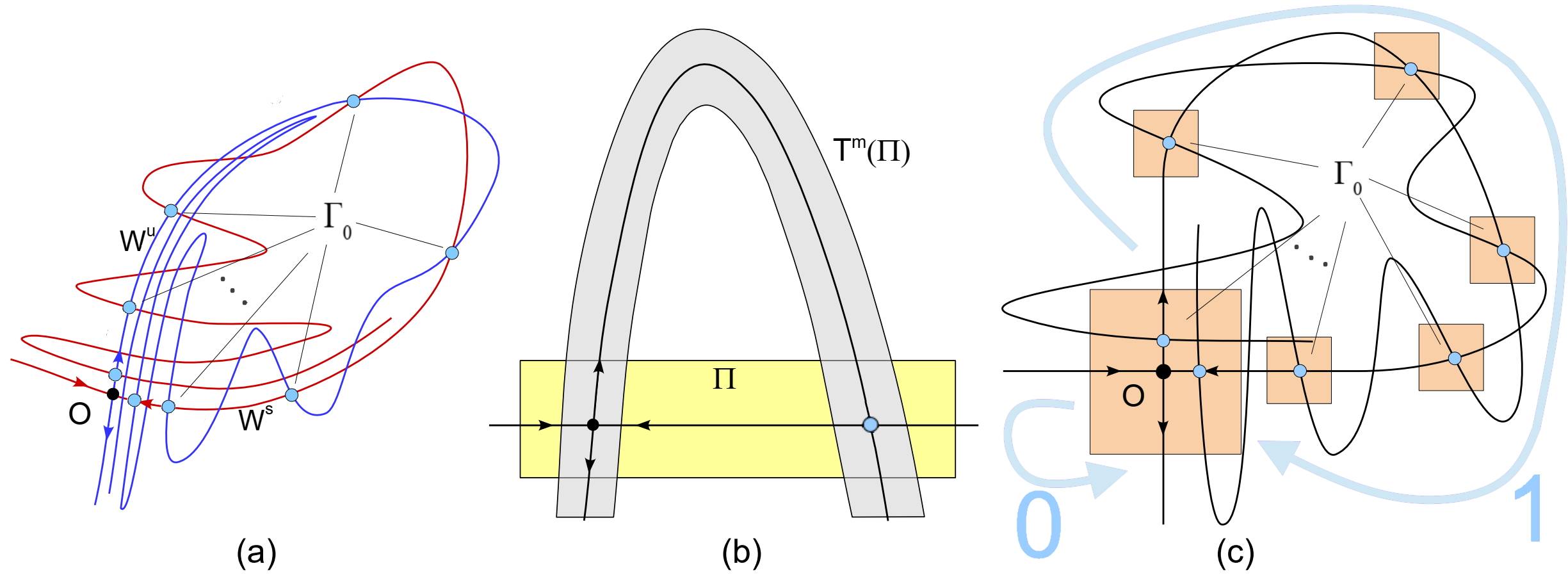}
\caption{\footnotesize (a) Poincar\'e homoclinic tangles. (b) Smale horseshoe
due to a transverse homoclinic intersection. (c) The Shilnikov method allows to account for all orbits in a small neighborhood (the orange boxes) of the transverse homoclinic orbit $\Gamma_0$. Each of these orbits is uniquely coded by a sequence of $0$s and $1$s: the symbol ``$0$'' corresponds to an iteration in a small neighborhood of the periodic orbit and ``$1$'' corresponds to an excursion along the homoclinic.}
\label{fig:Birk-1}
\end{figure}

One of the most fundamental statements of the theory of dynamical chaos concerning the dynamics near a {\em Poincar\'e homoclinic orbit}. This is an orbit which tends to the same saddle periodic orbit both in forward and backward time, i.e., it belongs to the intersections of the stable and unstable invariant manifolds $W^s$ and $W^u$ of the periodic orbit. A homoclinic orbit $\Gamma_0$ is called transverse when $W^s$ and $W^u$ cross transversely at the points of $\Gamma_0$, see Fig.~\ref{fig:Birk-1}a.

The possibility of the existence of transverse homoclinic orbits was established by H.~Poincar\'e in 1889\cite{Poin1}. This moment is universally accepted as the beginning of dynamical chaos history. Poincar\'e was awestruck by the complexity of the homoclinic tangles (see Fig.~\ref{fig:Birk-1}a) and stressed a special role of the transverse homoclinic as the universal mechanism of non-integrability. However a more or less detailed picture of the orbit behavior near the Poincar\'e homoclinic had remained unknown until mid 60s (the only partial result was obtained by D.~Birkhoff in 1934 who showed the existence of infinitely many periodic orbits for the special case where the map is two-dimensional and area-preserving\cite{Bir35}).

In 1965, S.~Smale \cite{Sm65} established the existence of a nontrivial hyperbolic set (the Smale horseshoe) near the transverse homoclinic (Fig.~\ref{fig:Birk-1}b). In 1967, Shilnikov closed the problem by providing a complete description of the structure of the set $N$ of orbits entirely lying in a small neighborhood of a transverse Poincar\'e homoclinic orbit. He proved \cite{Sh67Msb} that
\begin{itemize}
\item {\em the set $N$ is hyperbolic and is in one-to-one correspondence with the set of infinite sequences of two symbols (see Fig.~\ref{fig:Birk-1}c).}
\end{itemize}

We stress that Shilnikov was the first to pose and solve the problem of
a complete description of the dynamics near a transverse homoclinic orbit, thus demystifying the Poincare homoclinic tangles. He himself considered this result as fundamentally important and never tired of re-emphasizing that the Poincar\'e homoclinic orbit is ``an elementary building block'' of chaos, see his survey\cite{Msob03}.

One should also note that Smale's proof of the existence of the horseshoe was only done for a model case where the system near the saddle is assumed to be linear. However, this assumption cannot always be justified, e.g., it is not true for resonant saddles. Shilnikov's approach was free of this drawback.

To get rid of the linearization assumption, he had to overcome significant technical difficulties. Shilnikov had created a new mathematical technique for constructing solutions near saddle equilibrium states and periodic orbits, using the so-called ``method of boundary value problem'' (this method is described in detail in his original papers \cite{Sh67Msb, AfrSh73} and book \cite{book12}). He systematically applied this powerful method to the analysis of various homoclinic phenomena. He described homoclinic structures associated with invariant tori \cite{Sh68DAN} and, together with L.M.~Lerman, the dynamics near homoclinic orbits in infinitely-dimensional\cite{LerSh88} and general non-autonomous systems \cite{LerShChaos}. These works were far ahead of their time. We refer the reader to the recently published collection of selected Shilnikov's papers\cite{LPbook17}.

After these breakthrough results Shilnikov launched a new research program of the systematic study of plausible scenarios of transition to chaos. Since the early 70s, a vastly growing team of students and like-minded colleagues started getting centered around L.P.~Shilnikov. His first students were N.K.~Gavrilov, V.S.~Afraimovich, L.M.~Lerman, A.D.~Morozov, V.Z.~Grines, L.A.~Belyakov, and V.V.~Bykov, who themselves later became the world-known researchers. With time, the teams were further extended with V.I.~Lukyanov, Ya.L.~Umanskii, N.V.~Roschin, A.N.~Bautin, S.V.~Gonchenko, M.I.~Malkin, I.M.~Ovsyannikov, A.L.~Shilnikov,  V.S.~Biragov, D.V.~Turaev, M.V.~Shashkov, O.V.~Stenkin, I.V.~Belykh, V.S.~Gonchenko, and others, all from Nizhny Novgorod.  The rapidly growing activity of the Shilnikov scientific school ultimately led to the creation of a new branch in dynamical systems theory -- the theory of global bifurcations of high-dimensional systems.

The key challenge of this new theory was to identify and investigate the typical scenarios of the emergence of Poincare homoclinic orbits and, hence, chaotic dynamics. We have already discussed the transition to spiral chaos in Section \ref{spiral}. Below, we review the following major scenarios where Shilnikov made decisive contributions:
\begin{enumerate}
\item transition to chaos via a homoclinic tangency (a non-transverse Poincar\'e homoclinic orbit, see Section \ref{sec_3A});
\item the breakdown of an invariant torus (transition from quasi-periodic motions to chaos, see Section \ref{sec:tor});
\item the onset of the Lorenz attractor (see Section \ref{sec:lor}).
\end{enumerate}

\section{Homoclinic tangencies} \label{sec_3A}

It is well-known that the main property of chaotic dynamics is its structural instability (non-hyperbolicity): a typical chaotic system usually changes its behavior with arbitrarily small changes of parameters. A homoclinic tangency (a non-transverse intersection of invariant manifolds of a saddle periodic orbit) is a primary object in the theory of nonhyperbolic chaotic systems, see the book \cite{GS07}.

The systematic study of homoclinic tangencies was started by N.K. Gavrilov and L.P. Shilnikov in two papers\cite{GaS72,GaS73}. First of all, they classified different types of quadratic homoclinic tangencies and described the structure of {\em non-uniformly hyperbolic sets} entirely lying in a small neighborhood of the orbit of a homoclinic tangency. They also studied main bifurcations accompanying the splitting of the tangency in a generic one-parameter family and discovered that {\em stable periodic orbits, coexisting with hyperbolic sets}, can emerge at these bifurcations.
\begin{figure}[t!]
\includegraphics[width=.99\linewidth]{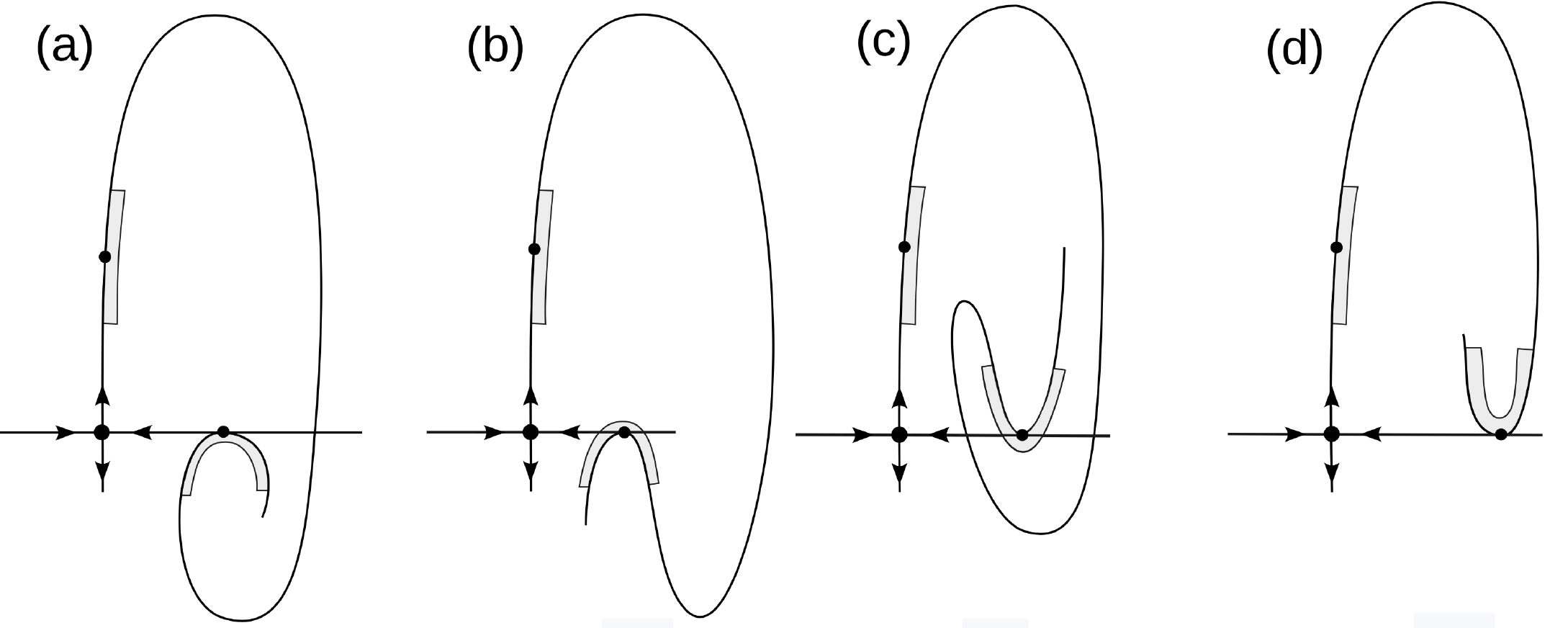}
\caption{Various types of homoclinic tangencies in 2D diffeomorphisms.}
\label{fig8}
\end{figure}

Let us illustrate the Gavrilov-Shilnikov results for the case of two-dimensional diffeomorphisms.  Let a saddle fixed point $O$ have multipliers $\lambda$, $\gamma$, such that $0 <|\lambda| <1 <|\gamma|$ and $\sigma \equiv |\lambda\gamma| < 1$. Let there exist a homoclinic orbit $\Gamma_0$ associated with a quadratic tangency of the stable and unstable manifolds $W^s(O)$ and $W^u(O)$. To be specific, assume $\lambda> 0$, $\gamma> 0$ (the cases with $\lambda$ and/or $\gamma$ negative were also considered in Refs.\cite{GaS72,GaS73}). Then, there are 4 different types of homoclinic tangencies, see Fig.~\ref{fig8}.

Systems with homoclinic tangencies {\em from below} (Figs.~\ref{fig8}a,b) belong to the {\it first class}. In this case, Gavrilov and Shilnikov proved that the set $N_0$ of all orbits that lie entirely in a small neighborhood $U$  of $O\cup \Gamma_0$ is trivial: $N_0 = \{O, \Gamma_0\}$. They also showed that
\begin{itemize}
\item systems with homoclinic tangencies of the first class can form the boundary between the systems with simple dynamics (the Morse-Smale systems) and systems with chaotic dynamics.
\end{itemize}
Thus, bifurcations of homoclinic tangencies of the first class (along with the bifurcation of a saddle-saddle with several homoclinic loops \cite{Sh69}) became the very first examples of what is known today as the {\it homoclinic $\Omega$-explosion}, when chaos (due to transverse homoclinics) emerge instantaneously, see more in Refs.\cite{NP76, PT85, StSh98, GSt11}.

In this connection, we also mention the review~\cite{Sh77} by Shilnikov, where he described key bifurcations that make Morse-Smale systems transition to chaos. One type of these bifurcations corresponds to an $\Omega$-explosion, as described above; saddle-node bifurcations that lead to the destruction of a two-dimensional invariant torus \cite{AfrSh74_6,AfrSh77} are another
example of the $\Omega$-explosion, see Section~\ref{sec:tor}. L.P. Shilnikov also described boundaries of the second kind that could be reached only through an infinite sequence (cascade) of bifurcations -- for example, he mentioned the period-doubling cascade that became famous after the work of M.~Feigenbaum\cite{F78}.

Systems with homoclinic tangencies {\em from above} (Figs.~\ref{fig8}c,d) belong either to the second or to the third class. A homoclinic tangency of the {\it second class} is shown in Fig.~\ref{fig8}c. In this case, Gavrilov and Shilnikov gave a {\em complete description} of $N_0$: it is a nontrivial non-uniformly hyperbolic set in one-to-one correspondence with the Bernoulli shift on three symbols 0,1,2 where the two homoclinic orbits $(...,0,...,0,1,0,...0,...)$ and $(...,0,...,0,2,0,...0,...)$ are glued together. They also showed that
\begin{itemize}
\item systems with homoclinic tangencies of the second class can form a boundary of the set of structurally stable systems with nontrivial uniformly-hyperbolic dynamics.
\end{itemize}

The partition of quadratic homoclinic tangencies into classes
was of great importance for the formation of the bifurcation theory of
chaotic dynamical systems. Depending on the tangency class, the global phase space transformations associated with the splitting of the tangency can lead to a drastic change in the system behavior -- from regular to chaotic, or from structurally stable to non-hyperbolic.

The tangency shown in Fig.~\ref{fig8}d is of the {\it third class}. Here, the set $N_0$ has a much more complicated structure than in the previous two cases. Gavrilov and Shilnikov described a large hyperbolic subset of $N_0$ and showed that
\begin{itemize}
\item the structure of the set $N_0$ depends essentially on the value of
$$
\theta = -\frac{\ln|\lambda|}{\ln|\gamma|},
$$
and, as a consequence, bifurcations of periodic orbits happen densely in
the set of systems with homoclinic tangencies of the third class, 
see also Refs.\cite{GS86,GS87}.
\end{itemize}

It was later shown by S.V.~Gonchenko and L.P.~Shilnikov \cite{G89,GS90,GS93} that the Gavrilov-Shilnikov parameter $\theta$ is an invariant (``modulus'') of the $\Omega$-equivalence (the topological equivalence on the non-wandering set) for systems with homoclinic tangencies of the third class (for other classes, $\theta$ is not an $\Omega$-modulus but, as follows from the results of J.~Palis on heteroclinic tangencies\cite{P78}, it is an invariant of topological equivalence). The further research revealed
\begin{itemize}
\item the existence of infinitely many independent $\Omega$-moduli for tangencies of the third class\cite{GST93a,GST99,G00}.
\end{itemize}

Informally, this means that one needs infinitely many independent parameters to describe the dynamics of any given system
with a homoclinic tangency of the third class. Therefore, for any
finite-parameters family of systems which has a homoclinic tangency at some parameter value, a complete description of the dynamics and bifurcations is not possible!

Another fundamental discovery in this theory was the paradoxical, at first glance, fact\cite{GST93a,GST99,GS07} that bifurcations of any system with a quadratic homoclinic tangency create {\em homoclinic tangencies of arbitrarily high orders}.
This means that the traditional logic of cancelling the degeneracy of a bifurcation by means of unfolding it within a parametric family of general position that depends on a sufficient (finite) number of parameters is not applicable to the study of non-hyperbolic chaos.

The following is a direct quote from L.P. Shilnikov\cite{shilnikov2000homoclinic,Msob03}:
\begin{itemize}
\item \textit{This discouraging result is even more important due to the denseness of the systems with homoclinic tangencies of the third class in the Newhouse regions, that is we have
the whole regions in the space of smooth dynamical systems for which a complete description of dynamics ... can never be achieved. When we realized all this, I remembered the words of E.A.~Leontovich concerning the discovery of chaos near a homoclinic loop to a saddle-focus: It just cannot be!}

\textit{In his report on the memoir of Poincar\'e, K.~Weierstrass wrote that the results of that paper eliminate many illusions of the theory of Hamiltonian systems. In essence, this was the starting point for the development of now qualitative methods which represent the essence of nonlinear dynamics. Today we see that the illusion of the possibility of a complete qualitative analysis of dynamical systems should also be abandoned. And in both cases, the cause of the crisis were Poincar\'e homoclinic curves.}
\end{itemize}

Since a complete description of non-hyperbolic chaos is beyond reach, one should concentrate on its most interesting properties. Such most general property is the coexistence of periodic orbits of different stability types. This fact, the emergence of the so-called ``windows of stability'' within chaos, is well known to anyone who performed numerical experiments with systems with strange attractors. The very first result about this belongs to Gavrilov and Shilnikov\cite{GaS73}:
\begin{itemize}
\item in the case $\sigma < 1$, for a generic one-parameter unfolding $X_{\mu}$ of a
two-dimensional map with a quadratic homoclinic tangency,
there exists a converging to $\mu=0$ infinite sequence of non-intersecting intervals $\delta_k$ of the values of the parameter $\mu$, such that $X_\mu$ has a stable orbit of period-$k$.
\end{itemize}
A few years later, S.~Newhouse showed\cite{N74,N79} that there are {\em open regions in the space of dynamical systems where systems with homoclinic tangencies are dense} and a generic system from the Newhouse region with $\sigma<1$ {\em has infinitely many stable periodic orbits whose closure contains a non-trivial hyperbolic set}.

Shilnikov knew well that the bifurcation of a homoclinic tangency is the most fundamental among bifurcations of chaotic systems. Therefore, the coexistence of hyperbolic sets with stable periodic orbits of sufficiently long periods is almost unavoidable for non-hyperbolic chaos. This, along with the work on the Lorenz model\cite{ABS77,ABS80,ABS82}, led him to the concept of a {\em quasiattractor}\cite{AfrSh83}, which we discuss in some
detail in Section \ref{sec:lor}.

Here, we would like to attract attention to another
theory rooted in this research. Replacing
$\sigma<1$ condition by $\sigma>1$ for a two-dimensional diffeomorphism $f$ is equivalent to replacing $f$ by its inverse, $f^{-1}$. Therefore, Gavrilov-Shilnikov's result on periodic sinks near a homoclinic tangency transforms to the problem on the existence of
periodic sources (unstable periodic orbits with all multipliers larger than $1$ in the absolute value). As discovered by S.V.~Gonchenko, L.P.~Shilnikov and D.V.~Turaev\cite{GST97}, there exist open regions in the space of two-dimensional diffeomorphisms, the so-called {\em absolute Newhouse regions}, where systems with homoclinic tangencies are dense with both $\sigma<1$ and $\sigma>1$. This implies that
\begin{itemize}
\item a generic diffeomorphism from an absolute Newhouse region has infinitely many periodic sinks, sources and saddles, and the closure of the set of sinks has a non-empty intersection (which contains
a non-trivial hyperbolic set) with the closure of the set of sources.
\end{itemize}

\begin{figure}[t!]
\includegraphics[width=.99\linewidth]{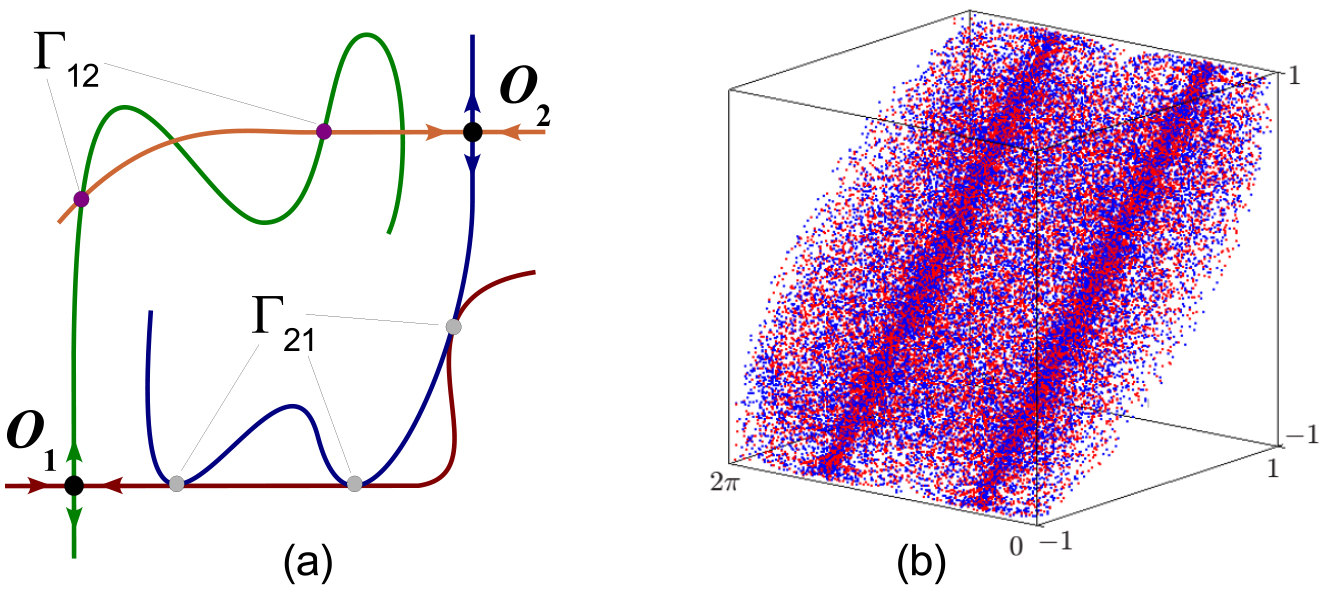}
\caption{{\footnotesize (a) A sketch of a 2D diffeomorphism $f_0$ with a non-transverse heteroclinic cycle containing two saddle fixed points $O_1$ and $O_2$ and two heteroclinic orbits $\Gamma_{12}$ and $\Gamma_{21}$ such that $0<J(O_1)<1<J(O_2)$, and $W^u(O_1)$ and $W^s(O_2)$ intersect transversely at the points of the orbit $\Gamma_{12}$, while $W^u(O_2)$ and $W^s(O_1)$ have a quadratic tangency at the points of the orbit $\Gamma_{21}$. (b) Celtic stone model~\cite{GGK13} produces the mixed dynamics: an attractor (red points) co-existing with a repeller (blue points) -- their intersection makes the picture to appear purple.}}
\label{fig15}
\end{figure}

The intervals of parameter values which belong to the absolute Newhouse regions exist\cite{GST97} in any generic one-parameter unfolding of a two-dimensional diffeomorphism with a non-transverse heteroclinic cycle which contains at least one saddle periodic orbit with the Jacobian (of the return map) greater than $1$ and one saddle periodic orbit with the Jacobian less than $1$, see Fig.~\ref{fig15}a.
Therefore, a non-hyperbolic chaotic map can fall into the absolute Newhouse domain whenever there are regions where the map expands areas and the regions where the map contracts areas and these regions are not dynamically separated.

Such situation is, in fact, quite common for reversible maps and for maps on closed surfaces. There are several competing definitions of what is ``an attractor of a dynamical system''. However, in any such definition, the attractor must contain all periodic sinks. Similarly, all periodic sources must belong to the repeller. Thus\cite{GST97}, 
\begin{itemize}
\item for systems from the absolute Newhouse domain the attractor and the repeller are inseparable in a persistent way.
\end{itemize}

For generalizations, e.g. to higher dimensions, see Refs.\cite{T96,GShSt02,GSS06,T10,T15}. Importantly, the robust intersection of the numerically observed attractor and repeller has been detected for
many examples\cite{GGK13,Kaz13,GGKT17,EN19,Kaz19, EN20,GGK20,Kaz20,Nek21}. This phenomenon was called
{\em mixed dynamics} in Ref.~\cite{GShSt02} and can, in addition to the classical conservative and dissipative chaos, be considered as a new type of chaotic dynamics\cite{GT17}.

\section{Mathematical theory of synchronization, and chaos through torus-breakdown} \label{sec:tor}

The classical theory of synchronization studies the effect of an external periodic forcing on a self-oscillatory system or the interaction of two self-oscillating systems. The goal of the theory is to identify synchronization regions in the parameter space corresponding to the existence of asymptotically stable periodic orbits and to to describe dynamical phenomena occurring on synchronization
boundaries.

The interest in this problem was originally motivated by applications, starting with the experimental studies by B.~Van der Pol and J.~Van der Mark \cite{VV27}, and by theoretical works by A.A.~Andronov and A.A.~Vitt \cite{AV30}. Mathematical models that appear in synchronization problems are multidimensional, and can exhibit nontrivial dynamics and chaos, as was, for the first time, discovered by M.L.~Cartwright and J.E.~Littlewood \cite{CL45,L57a,L57b}. Modern mathematical theory of synchronization started off in the fundamental works by V.S. Afraimovich and L.P. Shilnikov \cite{AfrSh74_4,AfrSh74_6,AfrSh77,AfrSh83} who proposed and developed new universal mechanisms of transitions from synchronized states to chaos. Prior to their work, it was mostly believed that the exit out of the synchronization region should necessarily lead to a quasi-periodic regime, i.e., to the appearance of a smooth two-dimensional invariant torus (by analogy with systems in a plane, where the disappearance of a homoclinic saddle-node results in the emergence of a stable limit cycle). Their papers destroyed the dogma by showing that the resonant torus could be non-smooth, and its breakdown could result in the onset of a chaotic invariant set \cite{AfrSh74_4,AfrSh74_6,AfrSh77,AfrSh83}.
   
When a dissipative self-oscillatory system is perturbed by an external signal of a small amplitude, then the initial limit cycle transforms into a two-dimensional asymptotically stable, smooth invariant torus $\tau$ in the extended phase space. The behavior of orbits on the torus can be studied using a two-parameter family $T_{\mu\omega}$ of Poincar\'e maps of some cross-section $S$ transverse to the torus. Here, parameters $\mu$ and $\omega$ are the amplitude and frequency of the external forcing, respectively.

\begin{figure}[h]
\includegraphics[width=.99\linewidth]{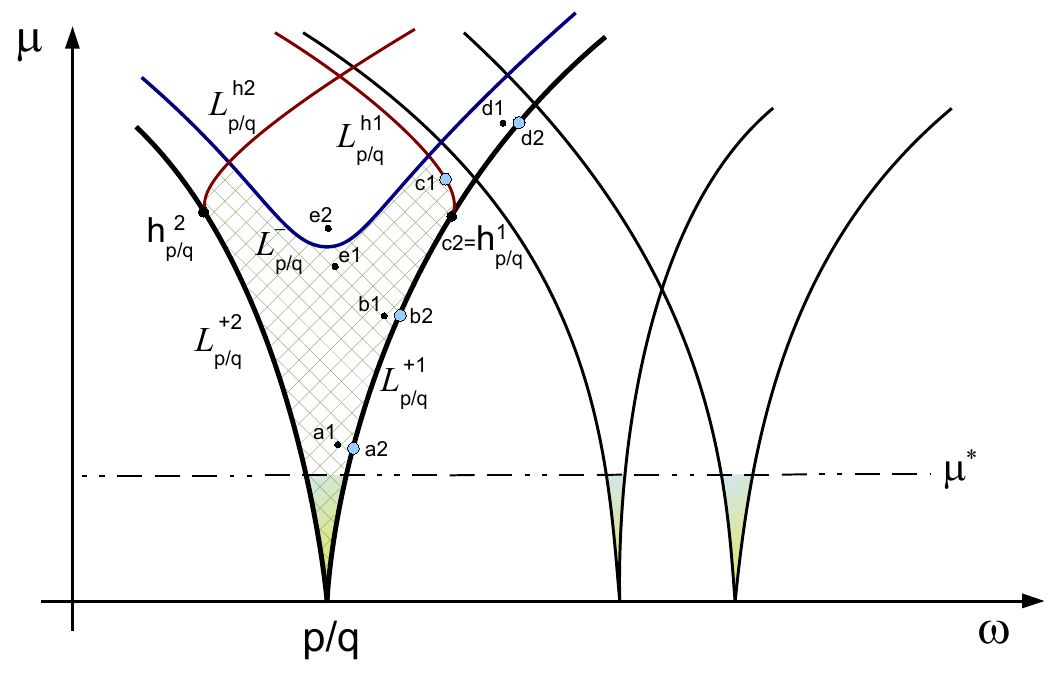}
\caption{Synchronization regions in the $(\omega,\mu)$-parameter plane.}
\label{fig9}
\end{figure}

The structure of the bifurcation set of the map $T_{\mu\omega}$ for small $\mu$ is the following\cite{Arn73,AAIS86}: for each rational $p/q$, there is the point $(\omega = p/q,\, \mu = 0)$ on the $\mu$-axis, which is a tip of the synchronization region $A_{p/q}$ in the $(\omega,\,\mu)$-parameter plane, as shown in Fig.~\ref{fig9}; it is also called the Arnold tongue. For sufficiently small $\mu$ ($\mu<\mu^*$ in Fig.~\ref{fig9}), different synchronization regions do not intersect, and for $(\mu,\omega) \in A_{p/q}$ the diffeomorphism $T_{\mu\omega}$ has a rational rotation number $p/q$. Thus, at least two periodic orbits exist on the invariant torus $\tau$. For simplicity, let there be exactly two such orbits. In the intersection with the cross-section $S$, they correspond to a pair of period $q$ orbits on the invariant curve $l_\mu =\tau\cap S$ -- one stable ($P_s$) and one saddle ($P_u$). The boundaries of the region $A_{p/q}$ consist of two bifurcation curves $L_{p/q}^{+1}$ and $L_{p/q}^{+2}$ which correspond to the existence of a saddle-node orbit of period $q$, formed as a result of the merger of $P_s$ and $P_u$, see Fig.~\ref{fig9}.

\begin{figure}[ht!]
\includegraphics[width=.99\linewidth]{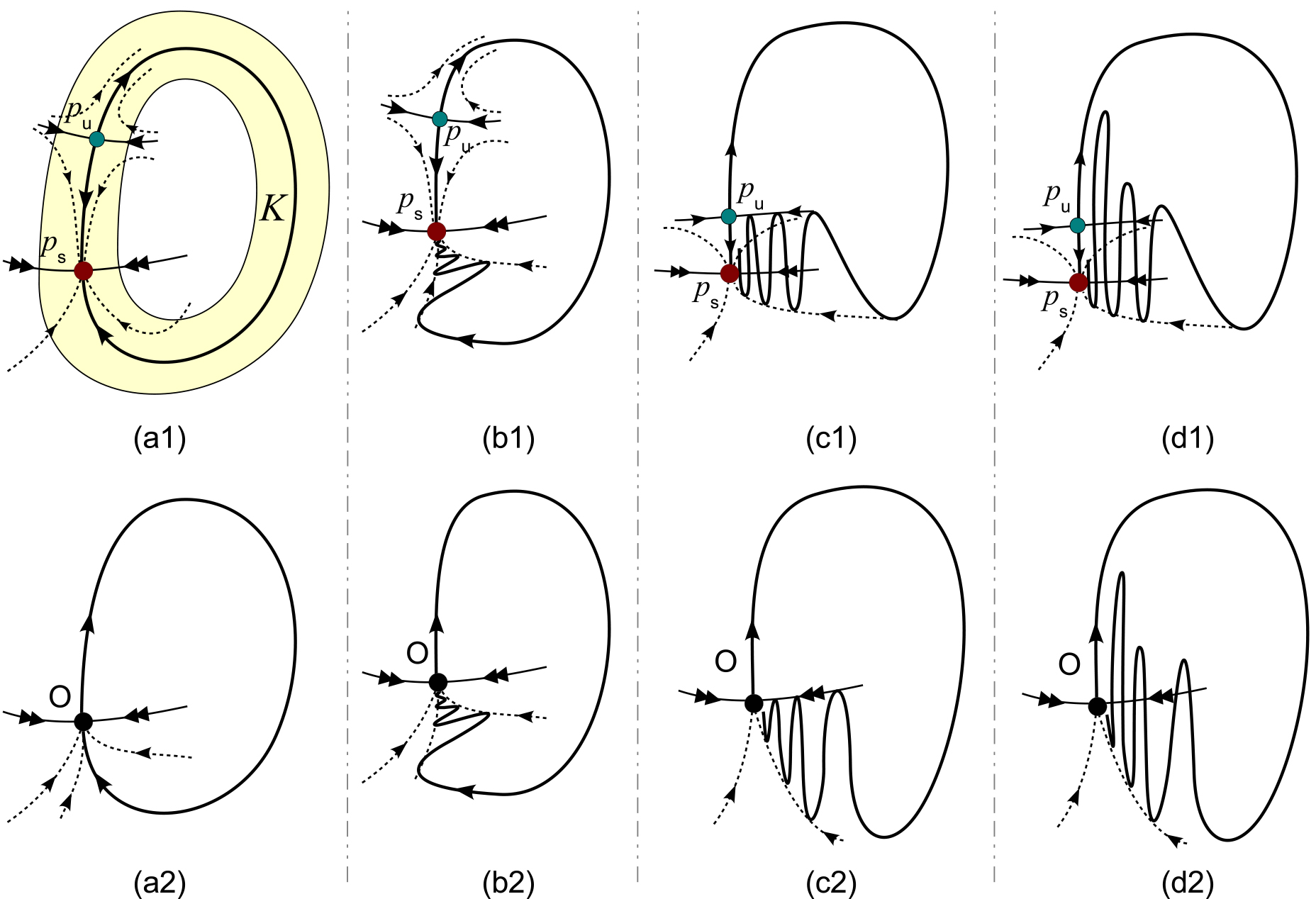}
\caption{The behavior of the unstable invariant manifold $W^u$ of the saddle point $P_u$ is shown in the figures in the upper row; the manifold $W^u(O)$ for the saddle-node $O$ is shown in the figures in the lower row. Parameter values are chosen from the corresponding regions in the bifurcation diagram of Fig.~\ref{fig9}; e.g., the upper row corresponds to the inside of the synchronization region and the lower row corresponds to its boundary.}
\label{fig10}
\end{figure}

To explain the Afraimovich-Shilnikov's discovery, we illustrate an evolution of the resonant circle $l_\mu$ in  Fig.~\ref{fig10} (for simplicity, we depict $P_s$ and $P_u$ as fixed points). Note that the invariant curve $l_\mu$ is the closure of the unstable manifold of $P_u $, i.e., $l_\mu = \overline{W^u(P_u)}$.

For small $\mu$, the curve $l_\mu$ is smooth: both the unstable separatrices of $P_u$ enter $P_s$ smoothly, see Fig.~\ref{fig10}(a1). However, as $\mu$ increases, the separatrices of $P_u$ can begin to oscillate and enter $P_s$ non-smoothly, as shown in Fig.~\ref{fig10}(b1). On the boundary of the synchronization region, the points $P_s$ and $P_u$ merge into a saddle-node $O$. If the invariant curve $l_\mu$ is smooth at this moment, see Fig.~\ref{fig10}(a2), the disappearance of the saddle-node does not break the invariant curve and dynamics remain quasi-periodic or periodic. However, if the unstable manifold of $O$ returns to it in a non-smooth way, see Fig.~\ref{fig10}(b2), then the exit from the synchronization region, can be accompanied by the emergence of chaotic dynamics. This is the main discovery of the pioneering works \cite{AfrSh74_4,AfrSh74_6,AfrSh77}, along with the very fact that the resonant torus can be non-smooth.

If the non-smoothness is strong enough, i.e., the so-called ``big loop condition'' is satisfied \cite{AfrSh74_4,AfrSh74_6,Sh86,SST04}, then chaos exists for all parameter values near the boundary outward of the synchronization region. When the non-smoothness is ``small'' and the big loop condition is not fulfilled, then intervals with chaotic and simple dynamics can alternate as the synchronization boundary is crossed over~\cite{NPT83,TSh86}.

Another related bifurcation, where the unstable manifold of a saddle-node has a transverse intersection with its strongly stable manifold, Fig.~\ref{fig10}(d2), was studied by V.I.~Lukyanov and L.P.~Shilnikov\cite{LukSh78}. This paper provided a theoretical foundation for the phenomenon of the alternation of regular and chaotic oscillations, which was later called the intermittency\cite{PM80}.

\begin{figure}[t!]
\includegraphics[width=.6\linewidth]{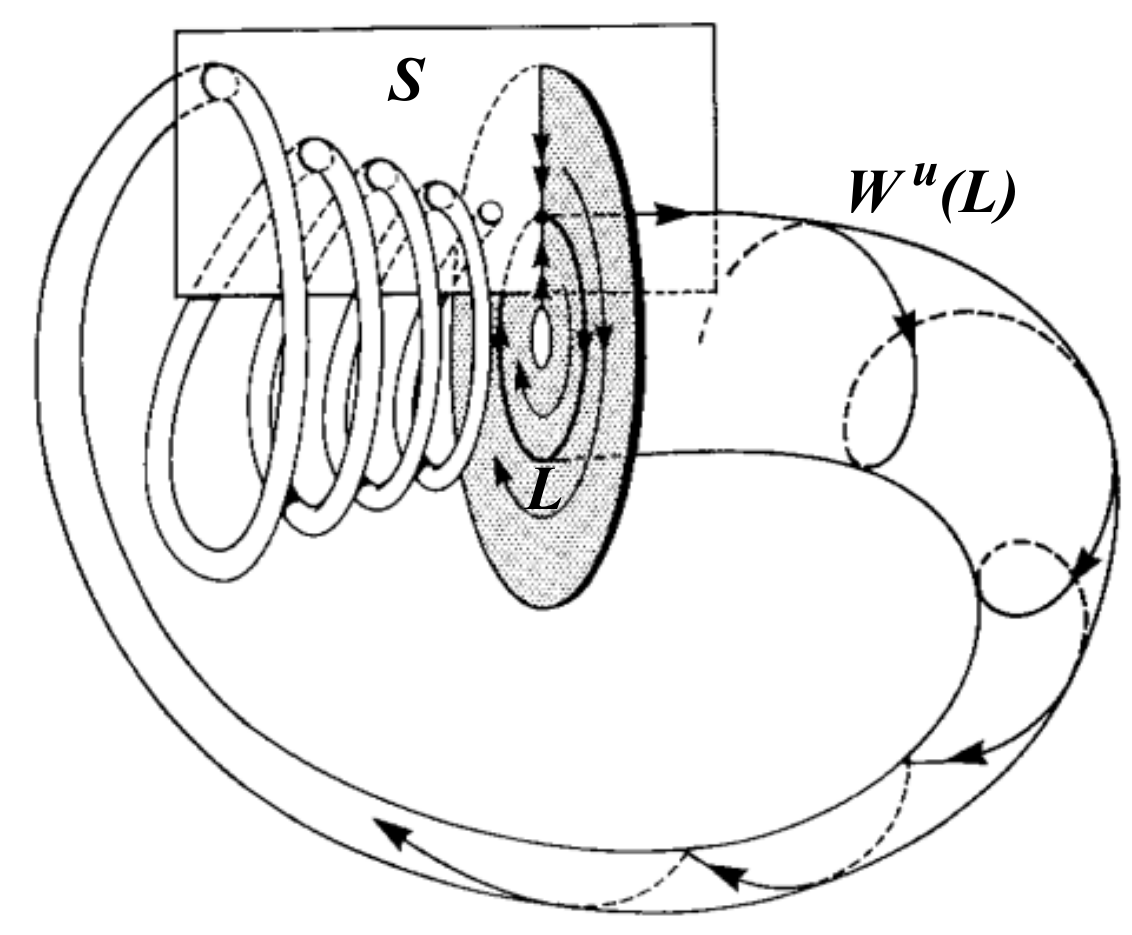}
\caption{Illustration to the mechanism of the blue-sky catastrophe, from Ref.~\cite{book12}. The squeezing unstable manifold $W^u(L)$ returns to the saddle-node periodic orbit $L$ from its node region so that the circles of its intersection with the cross-section $S$ shrink more and more with each subsequent iteration.}
\label{fig_bluescat}
\end{figure}

In the following papers\cite{TSh95,TSh00}, L.P.~Shilnikov and D.V.~Turaev presented a new type of homoclinic saddle-node bifurcations known as the ``blue sky catastrophe'', which is illustrated in Fig.~\ref{fig_bluescat}. This bifurcation produces a stable periodic orbit of an unbounded length. This was an example of a new stability boundary for periodic orbits, which was significantly different from the other eight known ones\cite{B49,Sh80a,book12}. Most probably, there are no more such boundaries of codimension-1, so the discovery of the blue-sky catastrophe completed the classical problem of the stability boundaries for periodic regimes. It was also shown\cite{TSh97}, that blue sky catastrophes offer new scenarios of the transition to chaos, e.g. some of them can lead to an instant formation of nontrivial hyperbolic attractors.

\begin{figure}[h!]
\includegraphics[width=.75\linewidth]{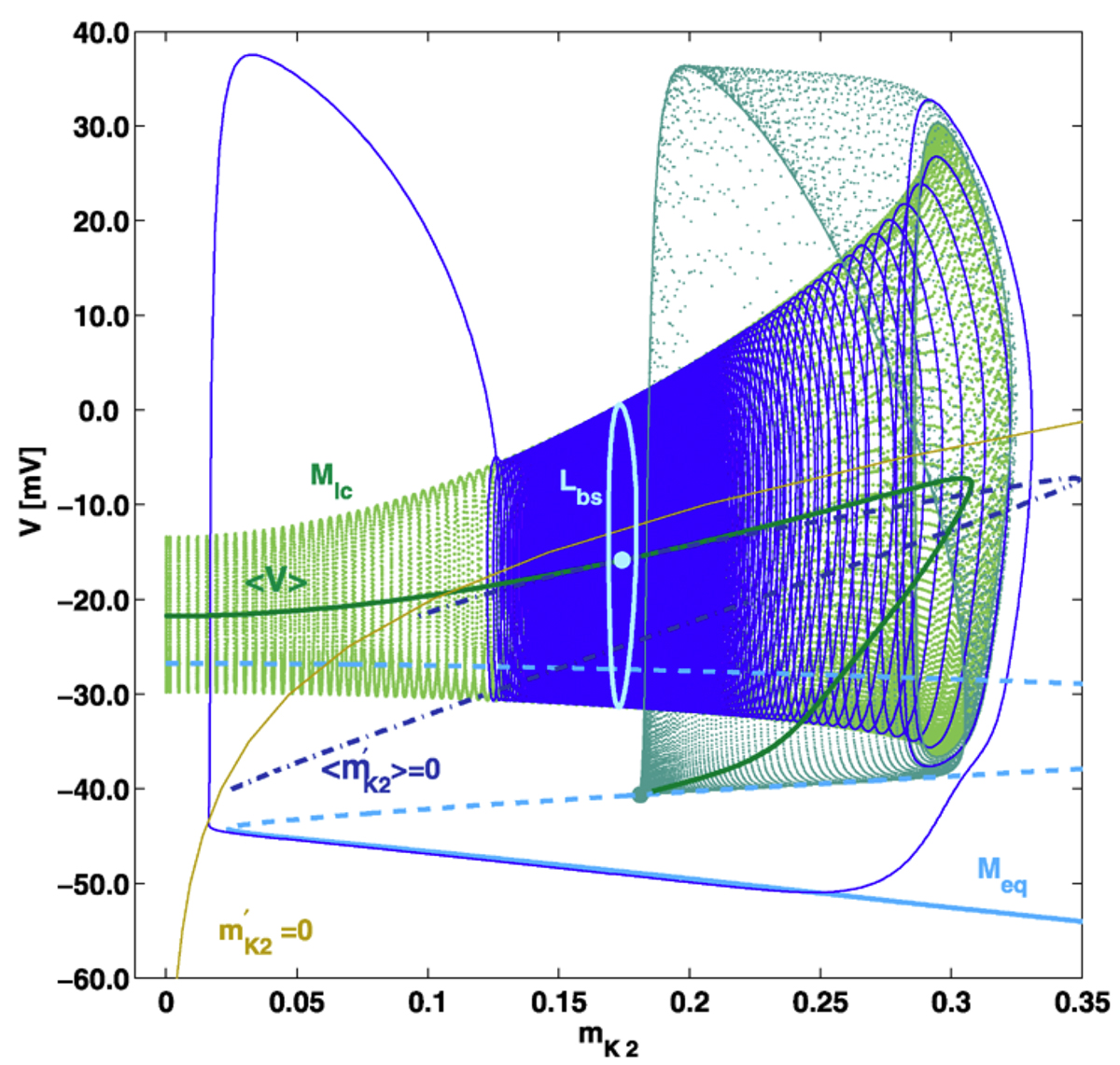}
\caption{Shown in the dark blue color is the long bursting orbit transitioning from spiking activity through blue sky catastrophe in the 3D phase space of the neural model as soon as the saddle-node limit cycle (in light blue) vanishes in the slow-motion manifold $M_{\rm LC}$; from Refs.~\cite{als2004,AL12}.}
\label{fig11}
\end{figure}

In further works, it was established that the blue sky catastrophe is a natural phenomenon in slow-fast systems \cite{SST05}, especially in those related to neuronal models \cite{cymbalyuk2005neuron, channell2007applications, shilnikov2014showcase}. This bifurcation is responsible for the transition from a regime of fast oscillations to the so-called bursting regimes, when a series of fast oscillations (bursts) follow one after another, alternated after slow quiescent transients\cite{SC05,SCC05}, as illustrated in Fig.~\ref{fig11}. It is also interesting that the homoclinic saddle-node bifurcation considered by Lukyanov and Shilnikov\cite{LukSh78}  turned out to be typical for many models of neurons, where it gives rise to complex bursting dynamics, as well as to the coexistence of fast spike and slow burst activity \cite{als2004,AL12}. The Afraimovich-Shilnikov torus bifurcations leading to complex dynamic activity have, too, happened to be universal phenomena for various low- and high-order models of neurons, see Refs.\cite{neural_torus,neural_chaos}  

\section{Lorenz attractor and beyond}
\label{sec:lor}
In the late 70s, deterministic chaos became a hot topic in nonlinear science. The question of to what degree the chaos discovered by mathematicians is relevant to natural sciences was widely discussed, and L.P.~Shilnikov took an active part in the conversation\cite{LeonSh70, Sh77, AfrShil1983, Sh84, Sh86,Sh85tb, Sh86}.

The turning point -- the de-facto proof that dynamical chaos is a fundamentally natural phenomenon -- was the discovery of a strange attractor in the Lorenz system~\cite{Lor63}. By that time, the only chaotic attractors known in the theory of dynamical systems were abstract examples of hyperbolic strange attractors. However, as these were not observed in applications, nonlinear scientists treated them merely as purely mathematical constructs. Once he learned about the Lorenz's work, Shilnikov realized right away the importance of the Lorenz model and that its study should lead to a new view on multidimensional dynamics. It was clear to Shilnikov that in his previous research on homoclinic bifurcations he already created all necessary machinery for studying Lorenz-type systems. He commissioned V.S.~Afraimovich and V.V.~Bykov to start the team work on the Lorenz attractor that resulted in a remarkable series of deep and influential publications.

\begin{figure}[t!]
\includegraphics[width=.99\linewidth]{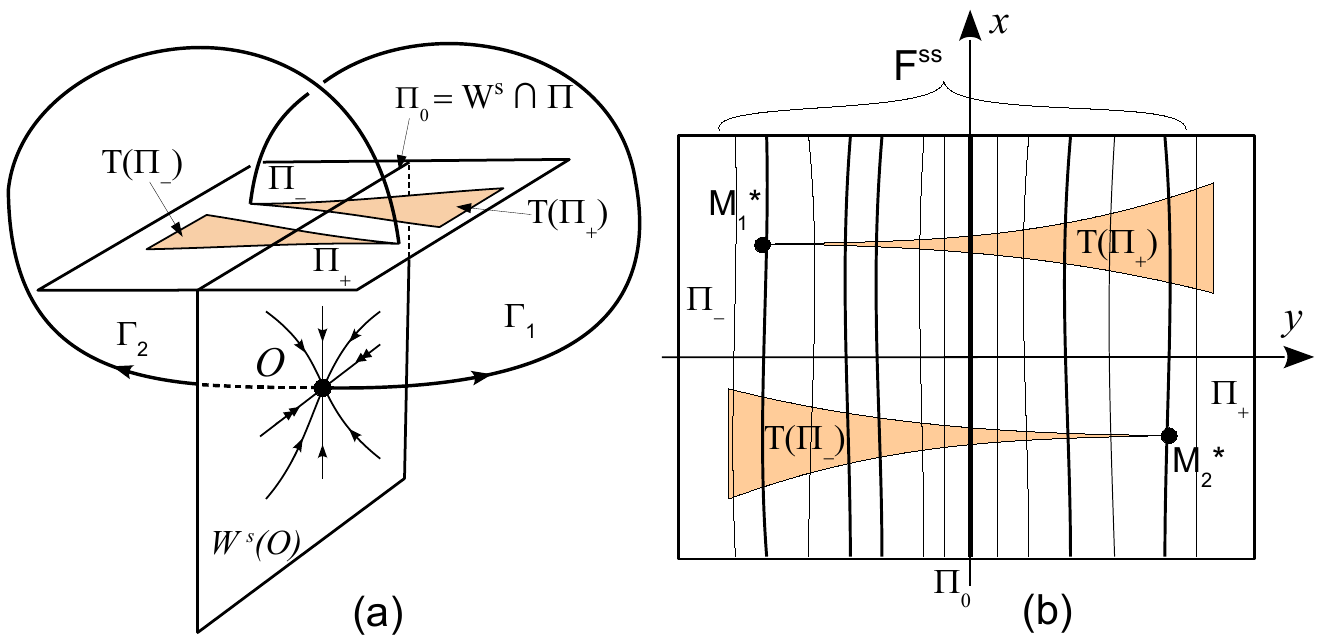}
\caption{
(a) Return map $T$ generated by trajectories of the Afraimovich-Bykov-Shilnikov geometric model. The map $T$ has a discontinuity line associated with the stable manifold of the saddle equilibrium $O$. It maps the cross-section $\Pi$ into a pair of triangles. (b) The strong-stable invariant foliation on $\Pi$ due to horizontal expansion and vertical contraction.}
\label{fig12}
\end{figure}

To rigorously study chaos in the Lorenz model, they introduced the so-called geometric model -- a two-dimensional piecewise smooth map whose derivatives satisfy certain specific inequalities. This is the Poincar\'e map on a cross-section $\Pi$ that
intersects the stable manifold of a saddle equilibrium state $O$; the map is discontinuous at the intersection line $\Pi_0$, see Fig. \ref{fig12}. Near the intersection line, the conditions on the derivatives are the same as for the Poincare map near a homoclinic loop to a saddle with positive saddle value and non-zero separatrix value\cite{Sh68a}, see Section \ref{sec:early}. Globally, these conditions have a form similar to the ``annulus principle'' proposed in L.P.~Shilnikov's earlier papers\cite{AfrSh74_6, AfrSh77} with V.S.~Afraimovich (on the breakdown of invariant tori, see Section \ref{sec:tor}). In essence, these are hyperbolicity conditions -- they express the contraction in the direction parallel to the discontinuity line $\Pi_0$ and expansion in the transverse direction. Importantly, the Afraimovich-Bykov-Shilnikov hyperbolicity conditions are explicitly verifiable. They were, first, numerically checked for the Lorenz model (at one value of parameters) by Ya.G. Sinai and E.B. Vul\cite{SinaiVul81} as well as by V.~Bykov and A.L.~Shilnikov \cite{BASh89} who determined  computationally the region  of the existence of the Lorenz attractor. At the classical Lorenz parameters ($\sigma=10$, $b=8/3$, $r=28$) the hyperbolicity conditions were verified by W.~Tucker\cite{Tucker99} who
used interval arithmetics to evaluate the numerical precision rigorously.

Using the geometric model, Afraimovich, Bykov and Shilnikov obtained a detailed description of the structure of the Lorenz attractor (see e.g. Theorems in Ref.\cite{ABS82}). They showed that the attractor is a (singularly) hyperbolic set and that in its neighborhood there exists a strongly contracting invariant foliation which makes the dynamics of the Poincare map effectively one-dimensional (i.e., it can be described by a piecewise expanding map of an interval\cite{Mal85,MalSaf21}). They also established that the attracting limit set consists of a two-dimensional transitive component where saddle periodic orbits are dense and, when the expansion is weak, it may also contain a one-dimensional component where the dynamics is described by a Markov chain. The latter case corresponds to the emergence of lacunas inside the attractor where the one-dimensional component resides. It was also shown in Ref.\cite{ABS82} that the Lorenz attractor is structurally unstable: the attractor itself persists but homoclinic loops of the saddle $O$ appear and disappear as parameters vary. Key bifurcations that lead to the formation of the Lorenz attractor, to the emergence of lacunas, and to the destruction of the attractor were also described\cite{ABS77,Sh80}.

The summary\cite{ABS77} of their findings was first published in 1977, which was followed by the detailed paper\cite{ABS82} containing the rigorous theory with complete proofs that came out 5 years (!) after its submission, due to bizarre circumstances, which had nothing to do with science. A groundbreaking survey~\cite{Sh80} of the results was prepared by~L.P. Shilnikov specifically to disseminate mathematical findings about the strange attractor and its bifurcations in the Lorenz model for the Russian-speaking nonlinear community.

Almost simultaneously and independently, a large number of publications on the Lorenz attractor appeared in the West \cite{lanford1977appendix,ruelle1976lorenz, guckenheimer1976strange, williams1979structure, yorke1979metastable, kaplan1979preturbulence, rand1978topological}, see also the collection of papers~\cite{SAsb}. Nevertheless, the Afraimovich-Bykov-Shilnikov theory remains the most complete and convenient for practical analysis of Lorenz-like attractors in various systems.

\begin{figure}[t!]
\includegraphics[width=.99\linewidth]{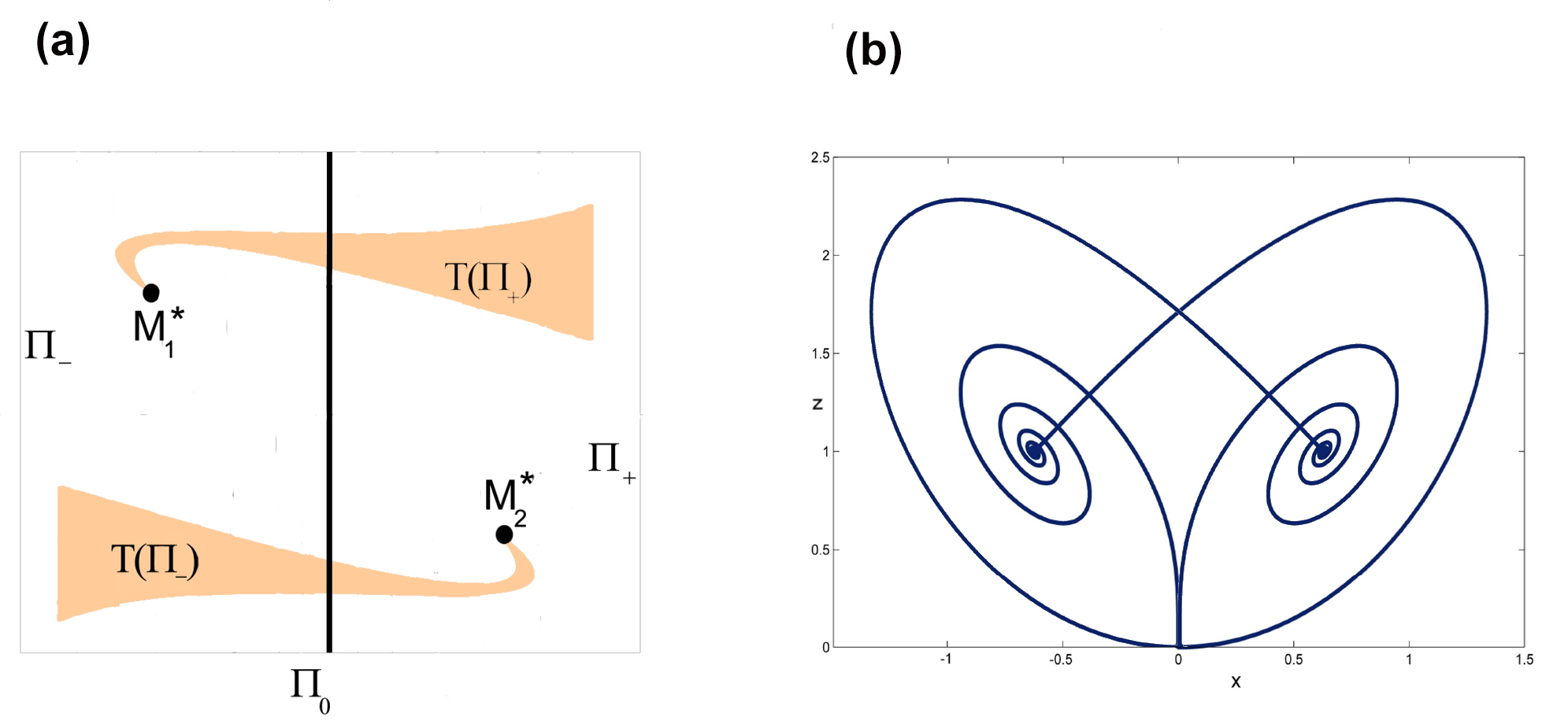}
\caption{(a) The ``hooks'' at the tips of $T(\Pi_+)$ and $T(\Pi_-)$ lead to homoclinic tangencies and, hence, to stable periodic orbits. (b) The heteroclinic cycle connecting a saddle and two symmetric saddle-foci (in the Shimizu-Morioka model\cite{ASh93}).}
\label{fig13}
\end{figure}

Shilnikov was very proud of these results, and the theme of the Lorenz attractor remained his priority until the end of his life\cite{TS08,GST09,BSS}. First and foremost, Shilnikov was interested in finding plausible bifurcation mechanisms leading to Lorenz-like attractors. He described a set of particular bifurcations of homoclinic loops\cite{Sh81} which warrant existence of the Lorenz attractor for an open set of parameter values nearby, without the need for
verifying the hyperbolicity conditions of the geometric model. These criteria
were subsequently used to examine the Lorenz attractor in the Shimizu-Morioka model\cite{ASh86,ASh89,ASh91,CTZ18} and in other systems\cite{robinson1989homoclinic,Rykhlik}. In the following paper\cite{SST93} the criteria were extended and new mechanisms of the onset
of the Lorenz attractor at local bifurcations were proposed.

The Lorenz attractor is perfectly chaotic: every orbit in it is unstable.
However, during the work on the Lorenz attractor papers\cite{ABS77,ABS82} it
had become clear to Shilnikov that this property does not always hold. In the Lorenz model itself, the change of parameters keeps the dynamics chaotic but, eventually, stable periodic orbits start emerging in the attractor, as was first reported in Ref.\cite{ABS80}. Shilnikov related this to the loss of the hyperbolicity due to the so-called ``hooks'' in the Poincar\'e map, see Fig.~\ref{fig13}a; the boundary in the parameter space of the Lorenz model where the hooks are formed was later found numerically in Ref.~ \cite{BASh89}.
It was noted in Ref.\cite{ABS80} that the formation of the hooks and the birth of stable periodic orbits is connected with the emergence of a symmetric heteroclinic cycle containing the saddle $O$ and a pair of saddle-foci, see Fig.~\ref{fig13}b. The study of the bifurcations of the heteroclinic cycles involving equilibria with different dimensions of the unstable manifolds was done by V.V.~Bykov\cite{Byk1,Byk2,Byk3} in his PhD thesis.

The fact\cite{ABS80} that by passing the ``Bykov point'' the strange attractor in the Lorenz model could loose its hyperbolicity property and gain stable periodic orbits served as an impetus for proposing the concept of a {\em quasiattractor} -- an attracting invariant set which, along with nontrivial hyperbolic subsets, can also contain long stable periodic orbits\cite{AfrShil1983}. These orbits emerged (via a homoclinic tangency, see Section \ref{sec_3A}) in almost any scenario of the development of chaos. Therefore, Shilnikov argued that the notion of the quasiattractors provides the most adequate mathematical description of dynamical chaos observed in most computational models.

\begin{figure}[t!]
\includegraphics[width=.99\linewidth]{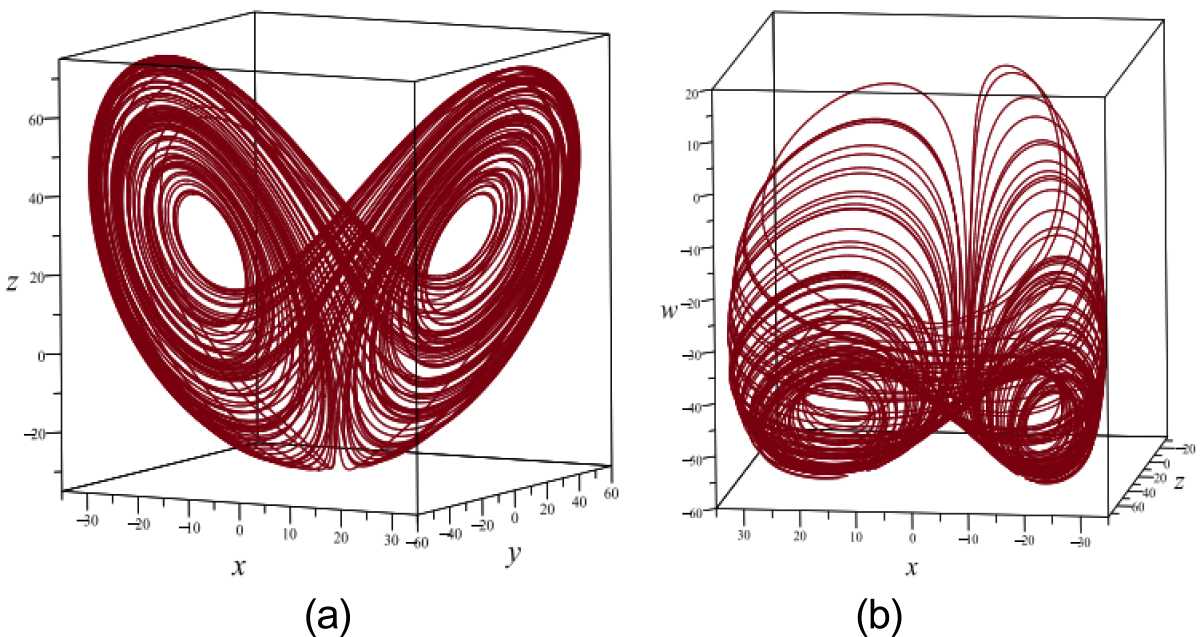}
\caption{(a) $(x,y,z)$-phase space and (b) the $(x,z,w)$-phase space projections of the strange attractor in the following 4D Lorenz-like system $\{\dot x = \sigma(y-x),\; \dot y = x(r-z)-y, \; \dot z = xy -bz + \mu y,\; \dot w = - b w - \mu z \}$ at $\sigma =10,\, b = 8/3,\, r = 25$ and $\mu  = 7$; from Ref.~\cite{GKT21}}
\label{fig14}
\end{figure}

The need to explore this paradigm stipulated the particular focus of many works by Shilnikov and his
research group on the study of global bifurcations that produce hyperbolic sets accompanied by stable periodic orbits. Thus, for multidimensional systems with homoclinic loops to a saddle-focus, I.M.~Ovsyannikov and L.P.~Shilnikov established\cite{ovsyannikov1986systems,ovsyannikov1991systems} the conditions for the appearance of stable periodic orbits near the homoclinic loop, as well as the specific criteria for their absence (both for the system itself and any system close to it).
The latter result led, subsequently, to the discovery of {\em wild spiral attractors}\cite{TS98} in  systems with the dimension four and higher. Such attractor is obtained from the Lorenz one when the saddle at the origin is replaced with a saddle-focus. Like the Lorenz attractor, the wild spiral attractor does not contain stable periodic orbits (due to a globalized version of the Ovsyannikov-Shilnikov criteria) and, moreover, all its orbits are unstable and remain such for any small perturbation of the system. Figure~\ref{fig14} shows the wild spiral attractor in a four-dimensional extension of the Lorenz system\cite{GKT21}.

Homoclinic tangencies occurring in the wild spiral attractor make a complete description of its structure and bifurcations impossible (see Section \ref{sec_3A}). Still, the main fact -- the robust
instability of every orbit in the attractor -- has
been established, based on its {\em pseudohyperbolicity} property. The concept and the theory of pseudohyperbolic strange attractors were proposed and developed\cite{TS98,TS08} by L.P.~Shilnikov and D.V.~Turaev.

Shilnikov proposed to apply this theory to the proof of chaoticity and to the further investigation of coupled Lorenz-like systems and periodically perturbed Lorenz-like systems\cite{TS08}. A particular example of pseudohyperbolic attractors that appear in this setting is given by {\em discrete Lorenz attractors}\cite{GOST05}. They look very similar to the classical Lorenz attractor but, as the name suggests, occur in systems with discrete time (i.e., diffeomorphisms). Such attractors have been detected in a broad class of applications\cite{GOST05,GGK13,GGOT13}. Universal bifurcation scenarios that explain why the discrete analogues of the Lorenz attractor are natural for
diffeomorphisms in dimension 3 and higher are described in
Ref.\cite{GGS12}.

\section{Content of the focus issue}

This Focus Issue contains 22 articles devoted to a variety of topical problems in the theory of dynamical systems and deterministic chaos.

The papers [\onlinecite{ASh21}, \onlinecite{MarG21}, \onlinecite{Kar21}, \onlinecite{SS21}] are focused on \textbf{spiral chaos}.

T.~Xing, K.~Pusuluri and A.L.~Shilnikov\cite{ASh21} reveal an intricate order of homoclinic bifurcations near the primary figure-8 connection of the Shilnikov saddle-focus in systems with central symmetry and reveal admissible shapes of the corresponding bifurcation curves in parameter space of such systems. They illustrate their theory using a newly developed symbolic toolbox to disclose the fine organization and self-similarity of bifurcation unfoldings in two model symmetric systems.

A.~Gonchenko, M.~Gonchenko, A.~Kozlov, and E.~Samylina \cite{MarG21} propose and investigate scenarios of the birth of discrete homoclinic attractors in three-dimensional non-orientable maps. These scenarios include the emergence of discrete spiral figure-8 attractors and Shilnikov attractors.

E.~Karatetskaia, A.~Shykhmamedov and A.~Kazakov \cite{Kar21} study in detail discrete Shilnikov homoclinic attractors for three-dimensional nonorientable Henon-like maps of the form $\bar x = y, \; \bar y = z, \; \bar z = Bx + Cy + Az - y^2$.

I.~Sataev and N.~Stankevich \cite{SS21} reveal scenarios of the hyperchaos formation in the modified Anishchenko-Astakhov generator. It is shown that these scenarios include bifurcation cascades leading to the emergence of discrete spiral Shilnikov attractors.\\

The following three papers~[\onlinecite{MS21}, \onlinecite{BBB21}, \onlinecite{GGKS21}] focus on \textbf{Lorenz-like attractors}.

M. Malkin and K. Safonov\cite{MS21} consider a two-parameter family of Lorenz-like maps of the interval. They have found the regions in the parameter plane where the topological entropy depends monotonically on the parameters, as well as the regions where the monotonicity is broken.

V.N. Belykh, N.V. Barabash and I.V. Belykh\cite{BBB21} study homoclinic bifurcations in a piecewise-smooth Lorenz-type system. They analytically construct the Poincar\'e return map and use it to establish the presence of sliding motions and, thereby, rigorously characterize sliding homoclinic bifurcations that destroy a chaotic Lorenz-type attractor.

S. Gonchenko, A. Gonchenko, A. Kazakov, and E. Samylina\cite{GGKS21} review geometrical and dynamical properties of discrete Lorenz-like attractors and propose new types of discrete pseudohyperbolic attractors.\\

The four papers [\onlinecite{KKS21}, \onlinecite{KK21}, \onlinecite{BGPY21}, \onlinecite{GM21}] are devoted to \textbf{hyperbolic attractors}.

V. Kruglov, P. Kuptsov and I.~Sataev\cite{KKS21} provide analytical and numerical evidence for the existence of a uniformly hyperbolic attractor of a Smale-Williams type  in a simple self-oscillating system with complex variables.

In the paper by P. Kuptsov and S.P. Kuznetsov\cite{KK21}, a slowly modulated self-oscillator system subject to a nonlinear delay is studied. It is shown that the system operates as two coupled hyperbolic chaotic subsystems and the transition to hyperbolic hyperchaos can occur.

M. Barinova, V. Grines, O. Pochinka, and B. Yu\cite{BGPY21} establish the existence of energy functions for three-dimensional diffeomorphisms with a hyperbolic attractor and a hyperbolic repeller.

V. Grines and D. Mints\cite{GM21} describe restrictions on possible   types of trivial and nontrivial hyperbolic basic sets of A-diffeomorphisms of surfaces.\\

The following four papers [\onlinecite{GS21}, \onlinecite{Kash21}, \onlinecite{DY21}, \onlinecite{MKBSO21}] discuss \textbf{synchronization and transition to chaos}.

I. Garashchuk and D. Sinelshchikov\cite{GS21} study the process of destruction of synchronous oscillations in a model of two interacting microbubble contrast agents exposed to an external ultrasound field.

S.A. Kashchenko\cite{Kash21} considers the effect of two types of unidirectional advective coupling in a ring chain of a large number of coupled Van der Pol equations. The local near-equilibrium dynamics are studied and the asymptotic behavior of the solutions is described.

Y. Deng and Y. Li\cite{DY21} propose a simple chaotic memristor-based circuit with an external stimulation and demonstrate its dynamical properties.

V. Munyaev, D. Khorkin, M. Bolotov, L. Smirnov, and G. Osipov \cite{MKBSO21} investigate chains of locally coupled identical pendulums with a constant torque. They show that chaos and hyperchaos appear as a result of changes in the individual properties of elements, and the properties of the entire ensemble under consideration.\\

The four papers~[\onlinecite{Tur21}, \onlinecite{Nek21}, \onlinecite{LT21}, \onlinecite{Mam21}] study {\em reversible and Hamiltonian} systems and the new type of dynamical chaos -- the \textbf{mixed dynamics}.

D. Turaev \cite{Tur21} gives criteria for non-conservative dynamics in reversible maps with transverse and non-transverse homoclinic orbits.

A. Emelianova and V. Nekorkin \cite{Nek21} describe the emergence of mixed dynamics in a system of two adaptively coupled phase oscillators under the action of a harmonic external force. They show that the mixed dynamics prevent the forced synchronization of the chaotic regime and, if an external force is applied to a reversible core, its fractal dimension decreases.

L. Lerman and K. Trifonov \cite{LT21} study homoclinic and heteroclinic bifurcations in reversible Hamiltonian systems with a saddle-center equilibrium and a saddle periodic orbit.

I. Bizyaev, S. Bolotin and I. Mamaev \cite{Mam21} investigate the dynamics of non-holonomic systems (the Chaplygin sleigh and the Suslov system) with periodically varying mass distribution.\\

The three papers [\onlinecite{Il21}, \onlinecite{MMP21}, \onlinecite{BGZ21}] are devoted to systems with simple dynamics, as well as the dynamics of foliations.

Yu. Iljashenko \cite{Il21} investigates the structure of the bifurcation diagrams of the families of vector fields in the plane. Among other results, he constructs a countable number of pairwise non-equivalent germs of two-parameter bifurcation diagrams.

Ya. Bazaikin, A. Galaev and N. Zhukova \cite{BGZ21} formulate the theory of the Cartan foliations and its relation to chaotic dynamics.

D. Malyshev, A. Morozov and O. Pochinka \cite{MMP21} propose a new approach to the classification of Morse-Smale diffeomorphisms on two-dimensional surfaces which allows for developing effective algorithms that can compute the topological class of a Morse-Smale diffeomorphism in polynomial time (as a function of the number of periodic points).

\begin{acknowledgments}
We would like to thank all contributors to this Focus Issue for their hard work and passion. Leonid Pavlovich Shilnikov would have enjoyed your papers.\\
A.K., S.G and D.T. acknowledge financial support from the Mathematics Science and Education Center ``Mathematics of Future Technologies,'' project no. 075-02-2021-1394.
\end{acknowledgments}

\section*{References}
%\bibliography{Shilnikov.bib}
\providecommand{\noopsort}[1]{}\providecommand{\singleletter}[1]{#1}%

\end{document}